\documentclass[reqno]{amsart}

\usepackage[usenames,dvipsnames]{pstricks}
\usepackage{epsfig}

\usepackage{color}
\date{\today}

\usepackage[latin1]{inputenc}
\usepackage[T1]{fontenc}
\usepackage{amsmath, amsthm}
\usepackage[english]{babel}
\usepackage{amssymb}
\usepackage{latexsym}
\usepackage{fancyhdr}
\usepackage{graphicx}
\usepackage[all]{xy}

\theoremstyle{defin}
\newtheorem{defin}{{\bf Definition}}[section]

\newtheorem{nota}[defin]{{\bf Remark}}
\newtheorem{exemplo}[defin]{{\bf Example}}

\newtheorem{teo}[defin]{{\bf Theorem}}

\newtheorem{corol}[defin]{{\bf Corollary}}
\newtheorem{lema}[defin]{{\bf Lemma}}

\newenvironment{dem}
{{\noindent{\bf Proof: }}\newline } {$\hfill\Box$\vspace{0.25cm}}

\newcommand{\ot}{\otimes}
\newcommand{\co}{\circ}

\begin{document}

\begin{center}
 {\huge{\bf  Iterated weak crossed products}}

\end{center}

\ \\

{\bf  J.M. Fern\'andez Vilaboa$^{1}$, R. Gonz\'{a}lez
Rodr\'{\i}guez$^{2}$ and A.B. Rodr\'{\i}guez Raposo$^{3}$}

\ \\
\hspace{-0,5cm}$1$ Departamento de \'Alxebra, Universidad de
Santiago de Compostela,  E-15771 Santiago de Compostela, Spain
(e-mail: josemanuel.fernandez@usc.es)
\ \\
\hspace{-0,5cm}$2$ Departamento de Matem\'{a}tica Aplicada II,
Universidad de Vigo, Campus Universitario Lagoas-Mar\-co\-sen\-de,
E-36310 Vigo, Spain (e-mail: rgon@dma.uvigo.es)
\\
\hspace{-0,5cm}$3$ Departamento de Matem\'aticas, Universidade da Coru\~{n}a,
Escuela Polit\'ecnica Superior, E-15403 Ferrol, Spain
(e-mail: abraposo@edu.xunta.es)
\ \\

\begin{center}
{\bf Abstract}
\end{center}
{\small In this paper we show how  iterate  weak crossed products with common monoid. More concretely, if $(A\ot V, \mu_{A\ot V})$ and $(A\ot W, \mu_{A\ot W})$ are weak crossed products, we find sufficient conditions 
to obtain a new weak crossed product $(A\ot V\ot W, \mu_{A\ot V\ot W})$.}

\vspace{0.5cm}



\section*{Introduction}

Let $A$ be a monoid and let $V$ be an object living in a strict
monoidal category $\mathcal C$ where every idempotent morphism
splits. In \cite{nmra4}  an associative product, called the weak crossed product of $A$ and
$V$, was defined on the tensor product $A\otimes V$  working with quadruples ${\Bbb A}_{V}=(A,V, \psi_{V}^{A}, \sigma_{V}^{A})$ where  $\psi_{V}^{A}:V\otimes A\rightarrow A\otimes V$
and $\sigma_{V}^{A}:V\otimes V\rightarrow A\otimes V$ are morphisms satisfying
some twisted-like and cocycle-like conditions. Associated to these
morphisms we define an idempotent morphism $\nabla_{A\otimes
V}:A\otimes V\rightarrow A\otimes V$ whose image,
denoted by $A\times V$, inherits the associative product from
$A\otimes V$. In order to define a unit for $A\times V$, and hence
to obtain a monoid structure in this object,  we complete the theory 
in \cite{mra-preunit}  using the notion of  preunit introduced by
Caenepeel and De Groot in \cite{Caen}. The theory presented in
\cite{nmra4} and \cite{mra-preunit} contains, as particular instances,  crossed
products where $\nabla_{A\otimes V}=id_{A\otimes
 V}$, for example the one defined by Brzezi\'nski in \cite{tb-crpr} or the notion of unified crossed product introduced by Agore and Militaru in \cite{AM2}, as well as crossed products where  $\nabla_{A\otimes V}\neq id_{A\otimes V}$ like, for example, the
weak smash product given by Caenepeel and De Groot in
\cite{Caen}, the notion of weak wreath products that we can
find in \cite{Street-WDL}, the weak crossed products for weak bialgebras
given in \cite{ana1} (see also \cite{mra-preunit}) and, as was proved in  \cite{mra-partial-unif}, the partial crossed products introduced by  Alves,
Batista, Dokuchaev and  Paques in \cite{partial}. Also,  B\"ohm
showed in \cite{bohm} that a monad in the weak version of the Lack
and Street's 2-category of monads in a 2-category is identical to a
crossed product system in the sense of \cite{nmra4}. Finally, weak crossed products appears in a natural way in the study of bilinear factorizations of algebras \cite{Gabi-Pepe1}, double crossed products of weak bialgebras \cite{Gabi-Pepe2}, and weak projections of weak Hopf algebras \cite{mra-proj}.

The purpose of this paper is to find an iteration process for weak crossed products with common monoid. Our main motivation comes from some interesting examples of this process that can be found in the recent literature. For example, in \cite{Pascual-Javi2}, Jara,  L\'opez, Panaite and Van Oystaeyen, motivated by the problem of defining a suitable representative for the product of spaces in noncommutative geometry, introduced the notion of  iterated
twisted tensor products of algebras. A good particular case of this iterated twisted tensor product  can be found in \cite{Maj}, where Majid constructed  an iterated sequence of double cross
products of certain bialgebras. On the other hand, in \cite{Pan-1}, Panaite proved that under
suitable conditions  a Brzezi\'nski crossed product may be iterated with a mirror version obtaining a new algebra structure. This construction contains as  examples the iterated twisted tensor product of algebras and the quasi-Hopf two-sided smash product. Finally, using the  2-category of weak distributive laws, B\"ohm describe in \cite{Bohm-iterated} a method of iterating Street's weak wreath product construction (see \cite{Street-WDL}). Note that in the two first examples of this paragraph the crossed products that we considered are cases where the associated idempotent is the identity. In that last one the associated idempotent it is not the identity. 

 An outline of the paper is as follows. Given two quadruples 
${\Bbb A}_{V}=(A,V, \psi_{V}^{A}, \sigma_{V}^{A})$ and ${\Bbb A}_{W}=(A,W, \psi_{W}^{A}, \sigma_{W}^{A})$ satisfying the suitable conditions that permit to obtain two weak crossed products  $(A\ot V, \mu_{A\ot V})$ and in $(A\ot W, \mu_{A\ot W})$, in the first section of this paper we introduce the notions of link and twisting morphism between ${\Bbb A}_{V}$ and ${\Bbb A}_{W}$, proving that, if they exist, it is possible to construct  a new quadruple ${\Bbb A}_{V\ot W}=(A, V\ot W, \psi_{V\ot W}^{A}, \sigma_{V\ot W}^{A}),$ satisfying the conditions that guarantee the existence of a new weak crossed product $(A\ot V\ot W, \mu_{A\ot V\ot W})$ called the iterated weak crossed product of $(A\ot V, \mu_{A\ot V})$ and $(A\ot W, \mu_{A\ot W})$. Also, we find the conditions under which there exists a preunit for $(A\ot V\ot W, \mu_{A\ot V\ot W})$. In the second section we discuss some examples involving wreath products, weak wreath products and the iteration process for Brzezi\'nski  crossed products proposed recently by D\u{a}u\c{s} and Panaite in \cite{Pan-2}. Finally in the last section, following the results proved in \cite{mra-preunit} we obtain a new characterisation of the iteration process.

Throughout this paper $\mathcal C$ denotes a strict  monoidal category with tensor product $\ot$, unit object $K$. There is no loss of generality in assuming that ${\mathcal C}$ is strict because by Theorem XI.5.3  of \cite{Christian} (this result implies the Mac Lane's coherence theorem) we know that every monoidal category is monoidally equivalent to a strict one. Then, we may work as if the constrains were all identities. We also assume that in ${\mathcal C}$ every idempotent morphism splits, i.e., for any morphism $q:M\rightarrow M$ such that $q\co q=q$ there exists an object $N$, called the image of $q$, and morphisms $i:N\rightarrow M$, $p:M\rightarrow N$ such that $q=i\co p$ and $p\co i=id_N$. The morphisms $p$ and $i$ will be called a factorization of $q$. Note that $Z$, $p$ and $i$ are unique up to isomorphism. The categories satisfying this property constitute a broad
class that includes, among others, the categories with epi-monic decomposition for
morphisms and categories with (co)equalizers.  Finally, given
objects $A$, $B$, $D$ and a morphism $f:B\rightarrow D$, we write
$A\ot f$ for $id_{A}\ot f$ and $f\ot A$ for $f\ot id_{A}$.

An monoid in ${\mathcal C}$ is a triple $A=(A, \eta_{A},
\mu_{A})$ where $A$ is an object in ${\mathcal C}$ and
 $\eta_{A}:K\rightarrow A$ (unit), $\mu_{A}:A\ot A
\rightarrow A$ (product) are morphisms in ${\mathcal C}$ such that
$\mu_{A}\co (A\ot \eta_{A})=id_{A}=\mu_{A}\co (\eta_{A}\ot A)$,
$\mu_{A}\co (A\ot \mu_{A})=\mu_{A}\co (\mu_{A}\ot A)$. Given two
monoids $A= (A, \eta_{A}, \mu_{A})$ and $B=(B, \eta_{B}, \mu_{B})$,
$f:A\rightarrow B$ is a monoid morphism if $\mu_{B}\co (f\ot
f)=f\co \mu_{A}$, $ f\co \eta_{A}= \eta_{B}$.

A comonoid in ${\mathcal C}$ is a triple ${D} = (D,
\varepsilon_{D}, \delta_{D})$ where $D$ is an object in ${\mathcal
C}$ and $\varepsilon_{D}: D\rightarrow K$ (counit),
$\delta_{D}:D\rightarrow D\ot D$ (coproduct) are morphisms in
${\mathcal C}$ such that $(\varepsilon_{D}\ot D)\co \delta_{D}=
id_{D}=(D\ot \varepsilon_{D})\co \delta_{D}$, $(\delta_{D}\ot D)\co
\delta_{D}=
 (D\ot \delta_{D})\co \delta_{D}.$ If ${D} = (D, \varepsilon_{D},
 \delta_{D})$ and
${ E} = (E, \varepsilon_{E}, \delta_{E})$ are comonoids,
$f:D\rightarrow E$ is a comonoid morphism if $(f\ot f)\co
\delta_{D} =\delta_{E}\co f$, $\varepsilon_{E}\co f
=\varepsilon_{D}.$

 Let  $A$ be a monoid. The pair
$(M,\varphi_{M})$ is a left $A$-module if $M$ is an object in
${\mathcal C}$ and $\varphi_{M}:A\ot M\rightarrow M$ is a morphism
in ${\mathcal C}$ satisfying $\varphi_{M}\circ( \eta_{A}\ot
M)=id_{M}$, $\varphi_{M}\circ (A\ot \varphi_{M})=\varphi_{M}\circ
(\mu_{A}\ot M)$. Given two left ${A}$-modules $(M,\varphi_{M})$ and
$(N,\varphi_{N})$, $f:M\rightarrow N$ is a morphism of left
${A}$-modules if $\varphi_{N}\circ (A\ot f)=f\circ \varphi_{M}$. In
a similar way we can define the notions of right $A$-module and
morphism of right $A$-modules. In this case we denote the left
action by $\phi_{M}$.

\section{Iterated weak crossed products}

In the first paragraphs of this section we resume some basic facts
about the general theory of weak crossed products. The complete details can be found in \cite{mra-preunit}.

Let $A$ be a monoid and $V$ be an object in
${\mathcal C}$. Suppose that there exists a morphism
$$\psi_{V}^{A}:V\ot A\rightarrow A\ot V$$  such that the following
equality holds
\begin{equation}\label{wmeas-wcp}
(\mu_A\ot V)\co (A\ot \psi_{V}^{A})\co (\psi_{V}^{A}\ot A) =
\psi_{V}^{A}\co (V\ot \mu_A).
\end{equation}
 As a consequence of (\ref{wmeas-wcp}), the morphism $\nabla_{A\ot V}:A\ot V\rightarrow
A\ot V$ defined by
\begin{equation}\label{idem-wcp}
\nabla_{A\ot V} = (\mu_A\ot V)\co(A\ot \psi_{V}^{A})\co (A\ot V\ot
\eta_A)
\end{equation}
is  idempotent. Moreover, $\nabla_{A\ot V}$ satisfies that
$$\nabla_{A\ot V}\co (\mu_A\ot V) = (\mu_A\ot V)\co
(A\ot \nabla_{A\ot V}),$$ that is, $\nabla_{A\ot V}$ is a left
$A$-module morphism (see Lemma 3.1 of \cite{mra-preunit}) for the
regular action  $\varphi_{A\ot V}=\mu_{A}\ot V$. With $A\times V$,
$i_{A\ot V}:A\times V\rightarrow A\ot V$ and $p_{A\ot V}:A\ot
V\rightarrow A\times V$ we denote the object, the injection and the
projection associated to the factorization of $\nabla_{A\ot V}$.
Finally, if $\psi_{V}^{A}$ satisfies (\ref{wmeas-wcp}), the following
identities hold
\begin{equation}\label{fi-nab}
(\mu_{A}\ot V)\co (A\ot \psi_{V}^{A})\co (\nabla_{A\ot V}\ot A)=
(\mu_{A}\ot V)\co (A\ot \psi_{V}^{A})=\nabla_{A\ot V}\co(\mu_{A}\ot
V)\co (A\ot \psi_{V}^{A}).
\end{equation}

From now on we consider quadruples ${\Bbb A}_{V}=(A, V,
\psi_{V}^{A}, \sigma_{V}^{A})$ where $A$ is a monoid, $V$ an
object, $\psi_{V}^{A}:V\ot A\rightarrow A\ot V$ a morphism
satisfiying (\ref{wmeas-wcp}) and $\sigma_{V}^{A}:V\ot V\rightarrow
A\ot V$  a morphism in ${\mathcal C}$.

We say that ${\Bbb A}_{V}=(A, V, \psi_{V}^{A}, \sigma_{V}^{A})$
satisfies the twisted condition if
\begin{equation}\label{twis-wcp}
(\mu_A\ot V)\co (A\ot \psi_{V}^{A})\co (\sigma_{V}^{A}\ot A) =
(\mu_A\ot V)\co (A\ot \sigma_{V}^{A})\co (\psi_{V}^{A}\ot V)\co
(V\ot \psi_{V}^{A})
\end{equation}
and   the  cocycle
condition holds if
\begin{equation}\label{cocy2-wcp}
(\mu_A\ot V)\co (A\ot \sigma_{V}^{A}) \co (\sigma_{V}^{A}\ot V) =
(\mu_A\ot V)\co (A\ot \sigma_{V}^{A})\co (\psi_{V}^{A}\ot V)\co
(V\ot\sigma_{V}^{A}).
\end{equation}

Note that, if ${\Bbb A}_{V}=(A, V, \psi_{V}^{A}, \sigma_{V}^{A})$
satisfies the twisted condition in Proposition 3.4 of
\cite{mra-preunit} we prove that the following equalities hold:
\begin{equation}\label{c1}
(\mu_A\otimes V)\circ (A\otimes \sigma_{V}^{A})\circ
(\psi_{V}^{A}\otimes V)\circ (V\otimes \nabla_{A\otimes V}) =
\nabla_{A\otimes V}\circ (\mu_A\otimes V)\circ (A\otimes
\sigma_{V}^{A})\circ (\psi_{V}^{A}\otimes V),
\end{equation}
\begin{equation}\label{aw}
\nabla_{A\otimes V}\circ (\mu_A\otimes V)\circ
(A\otimes\sigma_{V}^{A})\circ (\nabla_{A\otimes V}\otimes V) =
\nabla_{A\otimes V}\circ (\mu_A\otimes V)\circ
(A\otimes\sigma_{V}^{A}).
\end{equation}

Then, if $\nabla_{A\ot V}\co\sigma_{V}^{A}=\sigma_{V}^{A}$ we obtain
\begin{equation}\label{c11}
(\mu_A\otimes V)\circ (A\otimes \sigma_{V}^{A})\circ
(\psi_{V}^{A}\otimes V)\circ (V\otimes \nabla_{A\otimes V}) =
 (\mu_A\otimes V)\circ (A\otimes
\sigma_{V}^{A})\circ (\psi_{V}^{A}\otimes V),
\end{equation}
\begin{equation}\label{aw1}
 (\mu_A\otimes V)\circ
(A\otimes\sigma_{V}^{A})\circ (\nabla_{A\otimes V}\otimes V) =
(\mu_A\otimes V)\circ (A\otimes\sigma_{V}^{A}).
\end{equation}

By virtue of (\ref{twis-wcp}) and (\ref{cocy2-wcp}) we will consider
from now on, and without loss of generality, that
\begin{equation}
\label{idemp-sigma-inv} \nabla_{A\ot V}\co\sigma_{V}^{A} =
\sigma_{V}^{A}
\end{equation}
holds for all quadruples ${\Bbb A}_{V}=(A, V, \psi_{V}^{A},
\sigma_{V}^{A})$ {(see Proposition 3.7 of \cite{mra-preunit})}.

For ${\Bbb A}_{V}=(A, V, \psi_{V}^{A}, \sigma_{V}^{A})$ define the
product
\begin{equation}\label{prod-todo-wcp}
\mu_{A\ot  V} = (\mu_A\ot V)\co (\mu_A\ot \sigma_{V}^{A})\co (A\ot
\psi_{V}^{A}\ot V)
\end{equation}
and let $\mu_{A\times V}$ be the product
\begin{equation}
\label{prod-wcp} \mu_{A\times V} = p_{A\ot V}\co\mu_{A\ot V}\co
(i_{A\ot V}\ot i_{A\ot V}).
\end{equation}

If the twisted and the cocycle conditions hold, the product
$\mu_{A\ot V}$ is associative and normalized with respect to
$\nabla_{A\ot V}$ (i.e. $\nabla_{A\ot V}\co \mu_{A\ot V}=\mu_{A\ot
V}=\mu_{A\ot V}\co (\nabla_{A\ot V}\ot \nabla_{A\ot V}$)) and by the
definition of $\mu_{A\ot V}$ we have
\begin{equation}
\label{otra-prop} \mu_{A\ot V}\co (\nabla_{A\ot V}\ot A\ot
V)=\mu_{A\ot V}
\end{equation}
and therefore
\begin{equation}
\label{vieja-proof} \mu_{A\otimes V}\circ (A\otimes V\otimes
\nabla_{A\otimes V})=\mu_{A\otimes V}.
\end{equation}
 Due to the normality condition, $\mu_{A\times V}$ is
associative as well (Propostion 3.8 of \cite{mra-preunit}). Hence we
define:

\begin{defin}\label{wcp-def}{\rm
If ${\Bbb A}_{V}=(A, V, \psi_{V}^{A}, \sigma_{V}^{A})$  satisfies
(\ref{twis-wcp}) and (\ref{cocy2-wcp}) we say that $(A\ot V,
\mu_{A\ot V})$ is a weak crossed product.}
\end{defin}

The next natural question that arises is if it is possible to endow
$A\times V$ with a unit, and hence with a monoid structure. As
$A\times V$ is given as an image of an idempotent, it seems
reasonable to use the notion of  preunit introduced in \cite{Caen}
to obtain an unit. In our setting, if $A$ is a monoid, $V$ an
object in ${\mathcal C}$ and $m_{A\otimes V}$ is an associative
product defined in $A\otimes V$ a preunit $\nu:K\rightarrow A\otimes
V$ is a morphism satisfying
\begin{equation}
m_{A\otimes V}\circ (A\otimes V\otimes \nu)=m_{A\otimes V}\circ
(\nu\otimes A\otimes V)=m_{A\otimes V}\circ (A\otimes V\otimes
(m_{A\otimes V}\circ (\nu\otimes \nu))).
\end{equation}
Associated to a preunit we obtain an idempotent morphism
$$\nabla_{A\otimes
V}^{\nu}=m_{A\otimes V}\circ (A\otimes V\otimes \nu):A\otimes
V\rightarrow A\otimes V.$$ Take $A\times V$ the image of this
idempotent, $p_{A\otimes V}^{\nu}$ the projection and $i_{A\otimes
V}^{\nu}$ the injection. It is possible to endow $A\times V$ with a
monoid structure whose product is $$m_{A\times V} = p_{A\otimes
V}^{\nu}\circ m_{A\otimes V}\circ (i_{A\otimes V}^{\nu}\otimes
i_{A\otimes V}^{\nu})$$ and whose unit is $\eta_{A\times
V}=p_{A\otimes V}^{\nu}\circ \nu$ (see Proposition 2.5 of
\cite{mra-preunit}). If moreover, $m_{A\otimes V}$ is left
$A$-linear for the actions $\varphi_{A\otimes V}=\mu_{A}\otimes V$,
$\varphi_{A\otimes V\otimes A\otimes V }=\varphi_{A\otimes V}\otimes
A\otimes V$ and normalized with respect to $\nabla_{A\otimes
V}^{\nu}$,  the morphism
\begin{equation}
\label{beta-nu} \beta_{\nu}:A\rightarrow A\otimes V,\; \beta_{\nu} =
(\mu_A\otimes V)\circ (A\otimes \nu)
\end{equation}
is multiplicative and left $A$-linear for $\varphi_{A}=\mu_{A}$.

Although $\beta_{\nu}$ is not a monoid morphism, because $A\otimes
V$ is not a monoid, we have that $\beta_{\nu}\circ \eta_A = \nu$,
and thus the morphism $\bar{\beta_{\nu}} = p_{A\otimes
V}^{\nu}\circ\beta_{\nu}:A\rightarrow A\times V$ is a monoid
morphism.

In light of the considerations made in the last paragraphs, and
using the twisted and the cocycle conditions, in \cite{mra-preunit}
we characterize weak crossed products with a preunit, and moreover
we obtain a monoid structure on $A\times V$. These assertions are
a consequence of the following results proved in \cite{mra-preunit}.

\begin{teo}
\label{thm1-wcp} Let $A$ be a monoid, $V$ an object and
$m_{A\otimes V}:A\otimes V\otimes A\otimes V\rightarrow A\otimes V$
a morphism of left $A$-modules  for the actions $\varphi_{A\otimes
V}=\mu_{A}\otimes V$, $\varphi_{A\otimes V\otimes A\otimes V
}=\varphi_{A\otimes V}\otimes  A\otimes V$.

Then the following statements are equivalent:
\begin{itemize}
\item[(i)] The product $m_{A\otimes V}$ is associative with preunit
$\nu$ and normalized with respect to $\nabla_{A\otimes V}^{\nu}.$

\item[(ii)] There exist morphisms $\psi_{V}^{A}:V\otimes A\rightarrow A\otimes
V$, $\sigma_{V}^{A}:V\otimes V\rightarrow A\otimes V$ and
$\nu:k\rightarrow A\otimes V$ such that if $\mu_{A\otimes V}$ is the
product defined in (\ref{prod-todo-wcp}), the pair $(A\otimes V,
\mu_{A\otimes V})$ is a weak crossed product with $m_{A\otimes V} =
\mu_{A\otimes V}$ satisfying:
    \begin{equation}\label{pre1-wcp}
    (\mu_A\otimes V)\circ (A\otimes \sigma_{V}^{A})\circ
    (\psi_{V}^{A}\otimes V)\circ (V\otimes \nu) =
    \nabla_{A\otimes V}\circ
    (\eta_A\otimes V),
    \end{equation}
    \begin{equation}\label{pre2-wcp}
    (\mu_A\otimes V)\circ (A\otimes \sigma_{V}^{A})\circ
    (\nu\otimes V) = \nabla_{A\otimes V}\circ (\eta_A\otimes V),
    \end{equation}
    \begin{equation}\label{pre3-wcp}
(\mu_A\otimes V)\circ (A\otimes \psi_{V}^{A})\circ (\nu\otimes A) =
\beta_{\nu},
\end{equation}
\end{itemize}
where $\beta_{\nu}$ is the morphism defined in (\ref{beta-nu}). In
this case $\nu$ is a preunit for $\mu_{A\otimes V}$, the idempotent
morphism of the weak crossed product $\nabla_{A\otimes V}$ is the
idempotent $\nabla_{A\otimes V}^{\nu}$, and we say that the pair
$(A\otimes V, \mu_{A\otimes V})$ is a weak crossed product with
preunit $\nu$.
\end{teo}

\begin{nota}
\label{proof-resume} {\rm Note that in the proof of the previous
Theorem for $(i)\;\Rightarrow \;(ii)$ we define $\psi_{V}^{A}$ and
$\sigma_{V}^{A}$ as
\begin{equation}\label{fi-wcp}
\psi_{V}^{A} = m_{A\otimes V}\circ (\eta_A\otimes V\otimes
\beta_{\nu}),
\end{equation}
\begin{equation}\label{sigma-wcp}
\sigma_{V}^{A} = m_{A\otimes V}\circ (\eta_A\otimes V\otimes
\eta_A\otimes V).
\end{equation}
Also, by (\ref{pre3-wcp}), we have 
\begin{equation}\label{preunit-idemp}
\nabla_{A\ot V}\co \nu=\nu.
\end{equation}

}

\end{nota}

\begin{corol}\label{corol-wcp}
If $(A\otimes V, \mu_{A\otimes V})$ is a weak crossed product with
preunit $\nu$, then $A\times V$ is a monoid with the product
defined in (\ref{prod-wcp}) and unit $\eta_{A\times V}=p_{A\otimes
V}\circ\nu$.
\end{corol}

The aim of this section is to iterate weak crossed products with a
common monoid, that is weak crossed products induced by quadruples
of the form ${\Bbb A}_{V}=(A, V, \psi_{V}^{A}, \sigma_{V}^{A})$
where $A$ is fixed.

\begin{defin}\label{link}{\rm
Let ${\Bbb A}_{V}=(A, V, \psi_{V}^{A}, \sigma_{V}^{A})$ and ${\Bbb
A}_{W}=(A, W, \psi_{W}^{A}, \sigma_{W}^{A})$ be two quadruples. We
say that $$\Delta_{V\ot W}:V\ot W\rightarrow V\ot W$$ is a link
morphism between ${\Bbb A}_{V}$ and ${\Bbb A}_{W}$ if the following
conditions hold:
\begin{equation}
\label{falso-idemp} \psi_{V\ot W}^{A}=(A\ot \Delta_{V\ot W})\co \psi_{V\ot W}^{A},
\end{equation}
\begin{equation}
\label{falso-idemp2} \psi_{V\ot W}^{A}=\nabla_{A\ot V\ot W}\co(\psi_{V}^{A}\ot
W)\co (V\ot \psi_{W}^{A}),
\end{equation}
where 
$$
\psi_{V\ot W}^{A}=(\psi_{V}^{A}\ot W)\co (V\ot \psi_{W}^{A})\co
(\Delta_{V\ot W}\otimes A).
$$
and  $\nabla_{A\ot V\ot W}:A\ot V\ot W\rightarrow A\ot V\ot W$ is
the morphism defined by $$\nabla_{A\ot V\ot W}=(\mu_A\ot V\ot
W)\co(A\ot \psi_{V\ot W}^{A})\co (A\ot V\ot W\ot \eta_A).$$
 }
\end{defin}

\begin{lema}
\label{iden-ite}
 Let ${\Bbb A}_{V}=(A, V, \psi_{V}^{A},
\sigma_{V}^{A})$ and ${\Bbb A}_{W}=(A, W, \psi_{W}^{A},
\sigma_{W}^{A})$ be two quadruples. If there exists a link morphism
$\Delta_{V\ot W}:V\ot W\rightarrow V\ot W$ between them, the
morphism $\psi_{V\ot W}^{A}$ introduced in the previous definition
satisfies (\ref{wmeas-wcp}) and as a consequence $\nabla_{A\ot V\ot
W}$ is an idempotent morphism and the following identity holds:
\begin{equation}
\label{falso-idemp-link} \psi_{V\ot W}^{A}=\nabla_{A\ot V\ot
W}\co \psi_{V\ot W}^{A}.
\end{equation}

\end{lema}

\begin{dem} Using that $ \psi_{V}^{A}$, $ \psi_{W}^{A}$ satisfy
(\ref{wmeas-wcp}) and $\Delta_{V\ot W}$ satisfies
(\ref{falso-idemp}) we obtain

\begin{itemize}

\item[ ]$\hspace{0.38cm}(\mu_A\ot V\ot W)\co (A\ot \psi_{V\ot W}^{A})\co
(\psi_{V\ot W}^{A}\ot A)$

\item[ ]$=(\mu_A\ot V\ot W)\co (A\ot \psi_{V}^{A}\ot W)\co
(\psi_{V}^{A}\ot  \psi_{W}^{A})\co (V\ot \psi_{W}^{A}\ot A)\co
(\Delta_{V\ot W}\otimes A)$

\item[ ]$=(\psi_{V}^{A}\ot W)\co (V\ot \psi_{W}^{A})\co
(\Delta_{V\ot W}\ot \mu_{A}) $

\item[ ]$=\psi_{V\ot W}^{A}\co (V\ot W\ot \mu_{A})$

\end{itemize}
and then (\ref{wmeas-wcp}) holds for $\psi_{V\ot W}^{A}$. Finally, (\ref{falso-idemp-link}) follows directly from  (\ref{wmeas-wcp})  for $\psi_{V\ot W}^{A}$.

\end{dem}

\begin{defin}\label{wcp-def}{\rm
Let ${\Bbb A}_{V}=(A, V, \psi_{V}^{A}, \sigma_{V}^{A})$ and ${\Bbb
A}_{W}=(A, W, \psi_{W}^{A}, \sigma_{W}^{A})$ be two quadruples. We say that
$$\tau_{W}^{V}:W\ot V\rightarrow V\ot W$$ is a twisting morphism
between ${\Bbb A}_{V}$ and ${\Bbb A}_{W}$ if the following
conditions hold:

\begin{itemize}
\item[(i)] $(\psi_{V}^{A}\ot W)\co (V\ot
\psi_{W}^{A})\co (\tau_{W}^{V}\ot A)=(A\ot \tau_{W}^{V})\co
(\psi_{W}^{A}\ot V)\co (W\ot \psi_{V}^{A}).$

\item[ ]

\item[(ii)] $(\mu_{A}\ot V\ot W)\co (A\ot
\sigma_{V}^{A}\ot W)\co (\psi_{V}^{A}\ot \tau_{W}^{V})\co (V\ot
\sigma_{W}^{A}\ot V)\co (\tau_{W}^{V}\ot W\ot V)=$

\end{itemize}
$$(\mu_{A}\ot V\ot W)\co (A\ot \psi_{V}^{A}\ot
W)\co (A\ot V\ot \sigma_{W}^{A})\co (A\ot \tau_{W}^{V} \ot W)\co
(\psi_{W}^{A}\ot V\ot W)\co (V\ot \sigma_{V}^{A}\ot W)\co (W\ot V\ot
\tau_{W}^{V}).$$ }
\end{defin}

\begin{teo}
\label{prince}
 Let ${\Bbb A}_{V}=(A, V, \psi_{V}^{A},
\sigma_{V}^{A})$, ${\Bbb A}_{W}=(A, V, \psi_{W}^{A},
\sigma_{W}^{A})$ be two quadruples satisfying  (\ref{twis-wcp}) and
(\ref{cocy2-wcp}) with a link morphism $\Delta_{V\ot W}:V\ot
W\rightarrow V\ot W$ and with a twisting morphism $\tau_{W}^{V}:W\ot
V\rightarrow V\ot W$ between them. Then if we define $\sigma_{V\ot
W}^{A}:V\ot W\ot V\ot W\rightarrow A\ot V\ot W$ by \begin{equation}
\label{def-sigma} \sigma_{V\ot W}^{A}=(\mu_{A}\ot V\ot W)\co (A\ot
\psi_{V}^{A}\ot W)\co (\sigma_{V}^{A}\ot \sigma_{W}^{A})\co (V\ot
\tau_{W}^{V}\ot W)
\end{equation}
and it satisfies
\begin{equation}
\label{sigma1} \sigma_{V\ot W}^{A}=\sigma_{V\ot W}^{A}\co
(\Delta_{V\ot W}\ot V\ot W),
\end{equation}
\begin{equation}
\label{sigma2} \sigma_{V\ot W}^{A}=\sigma_{V\ot W}^{A}\co (V\ot W\ot
\Delta_{V\ot W}),
\end{equation}
\begin{equation}
\label{sigma3} \sigma_{V\ot W}^{A}= (A\ot \Delta_{V\ot W})\co
\sigma_{V\ot W}^{A},
\end{equation}
the quadruple
$${\Bbb A}_{V\ot W}=(A, V\ot W, \psi_{V\ot W}^{A}, \sigma_{V\ot W}^{A}),$$
where $\psi_{V\ot W}^{A}$ is the morphism defined in Lemma
\ref{iden-ite}, satisfies the equalities  (\ref{twis-wcp}),
(\ref{cocy2-wcp}) and (\ref{idemp-sigma-inv}). As a consequence,
$(A\ot V\ot W, \mu_{A\ot V\ot W})$ is a weak crossed product with
$$\mu_{A\ot V\ot W}=(\mu_A\ot V\ot W)\co (\mu_A\ot \sigma_{V\ot
W}^{A})\co (A\ot \psi_{V\ot W}^{A}\ot V).
$$
\end{teo}

\begin{dem} First we prove the twisted condition.

\begin{itemize}

\item[ ]$\hspace{0.38cm} (\mu_A\ot V)\co (A\ot \psi_{V\ot W}^{A})\co
(\sigma_{V\ot W}^{A}\ot A) $

\item[ ]$= (((\mu_{A}\ot V)\co (A\ot \psi_{V}^{A}))\ot W)\co (\sigma_{V}^{A}\ot ((\mu_{A}\ot W)
\co (A\ot \psi_{W}^{A})\co (\sigma_{W}^{A}\ot A))) \co (V\ot
\tau_{W}^{V}\ot W\ot A)  $
\item[ ]$\hspace{0.38cm} \co(\Delta_{V\ot W}\ot \Delta_{V\ot W}\ot A)$

\item[ ]$=((\mu_{A}\co (A\ot \sigma_{V}^{A})\co (\psi_{V}^{A}\ot V)\co (V\ot \psi_{V}^{A}))\ot W)\co
(V\ot V\ot (( \mu_{A}\ot W)\co (A\ot \sigma_{W}^{A})\co
(\psi_{W}^{A}\ot W)$
\item[ ]$\hspace{0.38cm}  \co (W\ot \psi_{W}^{A})))\co(V\ot \tau_{W}^{V}\ot
W\ot A)\co (\Delta_{V\ot W}\ot \Delta_{V\ot W}\ot A) $

\item[ ]$= (((\mu_{A }\ot V)\co (\mu_{A}\ot \sigma_{V}^{A})\ot
(A\ot \psi_{V}^{A}\ot V)\co (\psi_{V}^{A}\ot \psi_{V}^{A}))\ot W)
$
\item[ ]$\hspace{0.38cm}\co(V\ot (( \psi_{V}^{A}\ot \sigma_{W}^{A})\co
(V\ot \psi_{W}^{A}\ot W)\co (\tau_{W}^{V}\ot  \psi_{W}^{A}))) \co
(\Delta_{V\ot W}\ot \Delta_{V\ot W}\ot A) $

\item[ ]$=(((\mu_{A}\ot V)\co (A\ot ((\mu_{A}\ot V)\co (A\ot \sigma_{V}^{A})\co
(\psi_{V}^{A}\ot V)\co (V\ot \psi_{V}^{A}))))\ot W)  $
\item[ ]$\hspace{0.38cm} \co (A\ot V\ot V\ot \sigma_{W}^{A})\co (((\psi_{V}^{A}\ot \tau_{W}^{V})\co
(V\ot \psi_{W}^{A}\ot V)\co (V\ot W\ot \psi_{V}^{A}))\ot W) $
\item[ ]$\hspace{0.38cm} \co (V\ot W\ot V\ot \psi_{W}^{A})\co (\Delta_{V\ot W}\ot \Delta_{V\ot W}\ot A)  $

\item[ ]$= (\mu_A\ot V\ot W)\co (A\ot \sigma_{V\ot W}^{A})\co (\psi_{V\ot W}^{A}\ot V\ot W)\co
(V\ot W\ot \psi_{V\ot W}^{A}).$

\end{itemize}

In the previous calculus, the first equality follows by,
(\ref{sigma1}), (\ref{sigma2}), (\ref{sigma3}),  (\ref{wmeas-wcp})
for ${\Bbb A}_{V}$ and the associativity of $\mu_{A}$, the second
one follows by the twisted condition for ${\Bbb A}_{V}$ and ${\Bbb A}_{W}$,
the third one follows by (\ref{wmeas-wcp}) for ${\Bbb A}_{V}$ and the fourth
one follows  by (i) of Definition (\ref{wcp-def}) as well as the
associativity of $\mu_{A}$. Finally, in the last one we use the
twisted condition for ${\Bbb A}_{V}$.

The proof for the cocycle condition is the following:

\begin{itemize}

\item[ ]$\hspace{0.38cm}(\mu_A\ot V\otimes W)\co (A\ot \sigma_{V\otimes W}^{A}) \co
(\sigma_{V\otimes W}^{A}\ot V\otimes W) $

\item[ ]$= ((\mu_{A}\co (A\ot \mu_{A}))\ot V\ot W)\co (A\ot \mu_{A}\ot \psi_{V}^{A}\ot W)
\co (A\ot ((A\ot \sigma_{V}^{A})\co (\sigma_{V}^{A}\ot V))\ot A\ot
W) $
\item[ ]$\hspace{0.38cm}\co  (\psi_{V}^{A}\ot V\ot V\ot \sigma_{W}^{A})\co (V\ot
((\psi_{V}^{A}\ot \tau_{W}^{V})\co (V\ot \sigma_{W}^{A}\ot V)\co
(\tau_{W}^{V}\ot W\ot V))\ot W) $
\item[ ]$\hspace{0.38cm}\co  (\Delta_{V\ot W}\ot V\ot W \ot V\ot W)$

\item[ ]$=(\mu_{A}\ot V\ot W)\co (\mu_{A}\ot \psi_{V}^{A}\ot W)\co (A\ot \sigma_{V}^{A} \ot A\ot W)\co
(\psi_{V}^{A}\ot V\ot \sigma_{W}^{A}) $
\item[ ]$\hspace{0.38cm}\co  (V\ot [(\mu_{A}\ot V\ot W)\co (A\ot
\sigma_{V}^{A}\ot W)\co (\psi_{V}^{A}\ot \tau_{W}^{V})\co (V\ot
\sigma_{W}^{A}\ot V)\co (\tau_{W}^{V}\ot W\ot V)]\ot W)\co  $
\item[ ]$\hspace{0.38cm}\co  (\Delta_{V\ot W}\ot V\ot W \ot V\ot W)$

\item[ ]$=(\mu_{A}\ot V\ot W)\co (\mu_{A}\ot \psi_{V}^{A}\ot W)\co (A\ot \sigma_{V}^{A} \ot A\ot W)\co
(\psi_{V}^{A}\ot V\ot \sigma_{W}^{A}) $
\item[ ]$\hspace{0.38cm}\co (V\ot [(\mu_{A}\ot V\ot W)\co (A\ot \psi_{V}^{A}\ot
W)\co (A\ot V\ot \sigma_{W}^{A})\co (A\ot \tau_{W}^{V} \ot W)\co
(\psi_{W}^{A}\ot V\ot W) $
\item[ ]$\hspace{0.38cm}\co (V\ot \sigma_{V}^{A}\ot W)\co (W\ot V\ot
\tau_{W}^{V})]\ot W)\co   (\Delta_{V\ot W}\ot V\ot W \ot V\ot W)   $

\item[ ]$=(\mu_{A}\ot V\ot W)\co (\mu_{A}\ot \psi_{V}^{A}\ot W)\co (A\ot
((\mu_{A}\ot V)\co (A\ot \sigma_{V}^{A})\co (\psi_{V}^{A}\ot V)\co
(V\ot \psi_{V}^{A}))\ot \sigma_{W}^{A})  $
\item[ ]$\hspace{0.38cm}\co  (A\ot V\ot V\ot \sigma_{W}^{A}\ot W)\co
(\psi_{V}^{A}\ot \tau_{W}^{V} \ot W \ot W)\co (V\ot \psi_{W}^{A}\ot
V\ot W\ot W)  $
\item[ ]$\hspace{0.38cm} \co (V\ot W\ot ((\sigma_{V}^{A}\ot W)\co (V\ot \tau_{W}^{V}))\ot W)\co
  (\Delta_{V\ot W}\ot V\ot W \ot V\ot W)$

\item[ ]$=(\mu_{A}\ot V\ot W)\co (\mu_{A}\ot \psi_{V}^{A}\ot W)\co (A\ot
((\mu_{A}\ot V)\co (A\ot \psi_{V}^{A})\co (\sigma_{V}^{A}\ot V))\ot
\sigma_{W}^{A})  $
\item[ ]$\hspace{0.38cm}\co (A\ot V\ot V\ot \sigma_{W}^{A}\ot W)\co
(\psi_{V}^{A}\ot \tau_{W}^{V} \ot W \ot V)\co (V\ot \psi_{W}^{A}\ot
V\ot W\ot W)  $
\item[ ]$\hspace{0.38cm} \co (V\ot W\ot ((\sigma_{V}^{A}\ot W)\co (V\ot \tau_{W}^{V}))\ot W)
\co   (\Delta_{V\ot W}\ot V\ot W \ot V\ot W)$

\item[ ]$=(\mu_{A}\ot V\ot W)\co (\mu_{A}\ot \psi_{V}^{A}\ot W)\co (A\ot \sigma_{V}^{A}\ot
 ((\mu_{A}\ot W)\co (A\ot \sigma_{W}^{A})\co (\sigma_{W}^{A}\ot W)))  $
\item[ ]$\hspace{0.38cm} \co (((\psi_{V}^{A} \ot \tau_{W}^{V})\co (V\ot \psi_{W}^{A}\ot V)\co
(V\ot W\ot \sigma_{V}^{A}))\ot W\ot W)\co (V\ot W\ot V\ot
\tau_{W}^{V}\ot W)  $
\item[ ]$\hspace{0.38cm}\co (\Delta_{V\ot W}\ot V\ot W \ot V\ot W)$

\item[ ]$=(\mu_{A}\ot V\ot W)\co (\mu_{A}\ot \psi_{V}^{A}\ot W)\co (A\ot \sigma_{V}^{A}\ot
 ((\mu_{A}\ot W)\co (A\ot \sigma_{W}^{A})\co (\psi_{W}^{A}\ot W)\co (W\ot \sigma_{W}^{A})))  $
\item[ ]$\hspace{0.38cm} \co (((\psi_{V}^{A} \ot \tau_{W}^{V})\co (V\ot \psi_{W}^{A}\ot V)\co
(V\ot W\ot \sigma_{V}^{A}))\ot W\ot W)\co (V\ot W\ot V\ot
\tau_{W}^{V}\ot W)  $
\item[ ]$\hspace{0.38cm}\co (\Delta_{V\ot W}\ot V\ot W \ot V\ot W)$

\item[ ]$=(\mu_{A}\ot V\ot W)\co (A\ot \mu_{A}\ot V\ot W)\co (A\ot A\ot  \psi_{V}^{A}\ot W)
 $
\item[ ]$\hspace{0.38cm} \co  (A\ot((\mu_{A}\ot V)\co  (A\ot \psi_{V}^{A})\co (\sigma_{V}^{A}\ot
A))\ot \sigma_{W}^{A})\co(A\ot V\ot V\ot \psi_{W}^{A}\ot W)
$
\item[ ]$\hspace{0.38cm}\co  (\psi_{V}^{A}\ot \tau_{W}^{V}\ot \sigma_{W}^{A})\co  (V\ot
\psi_{W}^{A}\ot V\ot W\ot W)\co (V\ot W\ot ((\sigma_{V}^{A}\ot W)\co
(V\ot \tau_{W}^{V}))\ot W)$
\item[ ]$\hspace{0.38cm}\co (\Delta_{V\ot W}\ot V\ot W \ot V\ot W)$

\item[ ]$=(\mu_{A}\ot V\ot W)\co (A\ot \mu_{A}\ot V\ot W)\co (A\ot A\ot  \psi_{V}^{A}\ot W)
 $
\item[ ]$\hspace{0.38cm} \co  (A\ot((\mu_{A}\ot V)\co  (A\ot \sigma_{V}^{A})\co (\psi_{V}^{A}\ot
V)\co (V\ot \psi_{V}^{A}) )\ot \sigma_{W}^{A})\co(A\ot V\ot V\ot
\psi_{W}^{A}\ot W) $
\item[ ]$\hspace{0.38cm}\co (\psi_{V}^{A}\ot \tau_{W}^{V}\ot \sigma_{W}^{A})\co  (V\ot
\psi_{W}^{A}\ot V\ot W\ot W)\co (V\ot W\ot ((\sigma_{V}^{A}\ot W)\co
(V\ot \tau_{W}^{V}))\ot W)$
\item[ ]$\hspace{0.38cm}\co (\Delta_{V\ot W}\ot V\ot W \ot V\ot W)$

\item[ ]$=(\mu_{A}\ot V\ot W)\co (A\ot \mu_{A}\ot V\ot W)\co (A\ot A\ot  \psi_{V}^{A}\ot W)
 \co (\mu_{A}\co \sigma_{V}^{A}\ot \sigma_{W}^{A})$
\item[ ]$\hspace{0.38cm}\co (A\ot \psi_{V}^{A}\ot \tau_{W}^{V}\ot W)\co
( \psi_{V}^{A}\ot \psi_{W}^{A}\ot V\ot W)\co (V\ot \psi_{W}^{A}\ot
\psi_{V}^{A}\ot W)\co (V\ot W\ot \sigma_{V}^{A}\ot
\sigma_{W}^{A}) $
\item[ ]$\hspace{0.38cm} \co (V\ot W\ot V\ot \tau_{W}^{V}\ot W)\co
(\Delta_{V\ot W}\ot V\ot W \ot V\ot W)$

\item[ ]$=(\mu_A\ot V\otimes W)\co (A\ot \sigma_{V\otimes W}^{A})\co
(\psi_{V\otimes W}^{A}\ot V)\co
(V\otimes W\ot\sigma_{V\otimes W}^{A}). $

\end{itemize}

In this proof, the first equality follows by (\ref{sigma1}), the
associativity of $\mu_{A}$ and the twisted condition for ${\Bbb
A}_{V}$, the second one follows by the cocycle condition for ${\Bbb A}_{V}$,
(\ref{wmeas-wcp}) for ${\Bbb A}_{V}$ and the associativity of
$\mu_{A}$. In the third one  we used (ii) of Definition
(\ref{wcp-def}). The fourth one is a consequence of
(\ref{wmeas-wcp}) for ${\Bbb A}_{V}$ and the associativity of
$\mu_{A}$. The fifth one follows by the twisted condition for ${\Bbb
A}_{V}$ and the associativity of $\mu_{A}$. In the sixth one we used
(\ref{wmeas-wcp}) for ${\Bbb A}_{V}$ and the associativity of
$\mu_{A}$. The seventh one follows by the cocycle condition for
${\Bbb A}_{W}$ and in the eight one we applied (\ref{wmeas-wcp}) for
${\Bbb A}_{V}$ and the associativity of $\mu_{A}$ again. The ninth
one follows by the twisted condition for ${\Bbb A}_{V}$ and the
tenth one follows by  (i) of Definition (\ref{wcp-def}) and the
associativity of $\mu_{A}$. Finally, the last one was obtained using
(\ref{wmeas-wcp}) for ${\Bbb A}_{V}$ and ${\Bbb A}_{W}$.

The proof for the equality (\ref{idemp-sigma-inv}) is the following:

\begin{itemize}

\item[ ]$\hspace{0.38cm} \nabla_{A\ot V\ot W}\co \sigma_{V\ot W}^{A}$

\item[ ]$= (\mu_{A}\ot V\ot W)\co (A\ot \psi_{V}^{A}\ot W)\co (A\ot V\ot \psi_{W}^{A})\co
(((A\ot \Delta_{V\ot W})\co \sigma_{V\ot W}^{A})\ot \eta_{A}) $

\item [ ]$=(\mu_{A}\ot V\ot W)\co (A\ot \psi_{V}^{A}\ot W)\co (A\ot V\ot \psi_{W}^{A})\co
(\sigma_{V\ot W}^{A}\ot \eta_{A}) $

\item [ ]$=(\mu_{A}\ot V\ot W)\co (A\ot \psi_{V}^{A}\ot W)\co
 (\sigma_{V}^{A}\ot (\nabla_{A\ot W}\co  \sigma_{W}^{A}))\co
(V\ot \tau_{W}^{V}\ot W) $

\item [ ]$=\sigma_{V\ot W}^{A},$

\end{itemize}

where the first equality follows by definition, the second one by
(\ref{sigma2}), the third one by (\ref{wmeas-wcp}) for
$\psi_{V}^{A}$ and by the associativity of $\mu_{A}$. The last one relies on the properties of
$\sigma_{W}^{A}$, that is $\nabla_{A\ot W}\co \sigma_{W}^{A}=\sigma_{W}^{A}.$

\end{dem}

\begin{defin}\label{def-iter}{\rm
Let ${\Bbb A}_{V}=(A, V, \psi_{V}^{A}, \sigma_{V}^{A})$, ${\Bbb
A}_{W}=(A, W, \psi_{W}^{A}, \sigma_{W}^{A})$ be two quadruples
satisfying  (\ref{twis-wcp}) and (\ref{cocy2-wcp}) with a link
morphism $\Delta_{V\ot W}:V\ot W\rightarrow V\ot W$ and with a
twisting morphism $\tau_{W}^{V}:W\ot V\rightarrow V\ot W$ between
them. Let $(A\ot V, \mu_{A\ot V})$ and $(A\ot W, \mu_{A\ot W})$ be the
weak crossed products associated to ${\Bbb A}_{V}$ and ${\Bbb
A}_{W}$ and suppose that the morphism $\sigma_{V\ot
W}^{A}$ defined in (\ref{def-sigma}) satisfies (\ref{sigma1}),
(\ref{sigma2}) and (\ref{sigma3}). The weak crossed product $(A\ot
V\ot W, \mu_{A\ot V\ot W})$ defined in the previous theorem will be
called the iterated weak crossed product of $(A\ot V, \mu_{A\ot V})$
and $(A\ot W, \mu_{A\ot W})$. }

\end{defin}

In the following theorem we introduce the conditions that implies the existence of a preunit for the iterated weak crossed product defined previously.

\begin{teo}
\label{Teo-preunit}
Let ${\Bbb A}_{V}=(A, V, \psi_{V}^{A}, \sigma_{V}^{A})$, ${\Bbb
A}_{W}=(A, W, \psi_{W}^{A}, \sigma_{W}^{A})$ be two quadruples
satisfying  (\ref{twis-wcp}) and (\ref{cocy2-wcp}) with a link
morphism $\Delta_{V\ot W}:V\ot W\rightarrow V\ot W$ and with a
twisting morphism $\tau_{W}^{V}:W\ot V\rightarrow V\ot W$ between
them. Let $(A\ot V, \mu_{A\ot V})$ and $(A\ot W, \mu_{A\ot W})$ be the
weak crossed products associated to ${\Bbb A}_{V}$ and ${\Bbb
A}_{W}$  and suppose that $\nu_{V}:K\rightarrow A\ot V$
and $\nu_{W}:K\rightarrow A\ot W$ are preunits for $\mu_{A\ot V}$
and $\mu_{A\ot W}$. If  the morphism $\sigma_{V\otimes
W}^{A}$ defined in (\ref{def-sigma}) satisfies (\ref{sigma1}),
(\ref{sigma2}), (\ref{sigma3}) and the following equalities hold
\begin{equation}
\label{pre-1} (\mu_{A}\ot V\ot W)\co (A\ot \sigma_{V}^{A}\ot W)\co
(\psi_{V}^{A}\ot \tau_{W}^{V})\co (V\ot \psi_{W}^{A}\ot V)\co
(\Delta_{V\ot W}\ot \nu_{V})
\end{equation}
$$=\nabla_{A\ot V\ot W}\co (\eta_{A}\ot V\ot W),$$
\begin{equation}
\label{pre-2} (\mu_{A}\ot V\ot W)\co (A\ot \psi_{V}^{A}\ot W)\co
(A\ot V\ot \sigma_{W}^{A})\co (A\ot \tau_{W}^{V}\ot W)\co
(\nu_{W}\ot V\ot W)
\end{equation}
$$ =\nabla_{A\ot V\ot W}\co (\eta_{A}\ot V\ot W),$$
the iterated weak crossed product of $(A\ot V, \mu_{A\ot V})$ and
$(A\ot W, \mu_{A\ot W})$ has a preunit defined by
\begin{equation}
\label{iterated-preunit} \nu_{V\ot W}=\nabla_{A\ot V\ot W}\co
(\mu_{A}\ot V\ot W)\co (A\ot\psi_{V}^{A}\ot W)\co (\nu_{V}\ot
\nu_{W}).
\end{equation}
\end{teo}

\begin{dem} Note that to prove that $\nu_{V\ot W}$ is a preunit we
need to show that the equalities (\ref{pre1-wcp}), (\ref{pre2-wcp})
and (\ref{pre3-wcp}) hold for the quadruple ${\Bbb A}_{V\ot W}=(A,
V\ot W, \psi_{V\ot W}^{A}, \sigma_{V\ot W}^{A}).$

In this setting, the equality (\ref{pre3-wcp}) holds because:

\begin{itemize}

\item[ ]$\hspace{0.38cm}(\mu_{A}\ot V\ot W)\co (A\ot \psi_{V\ot W}^{A})\co (\nu_{V\ot W}\ot A)$

\item[ ]$= \nabla_{A\ot V\ot W}\co (\mu_{A}\ot V\ot W)\co (A\ot \mu_{A}\ot V\ot W)\co
(A\ot A\ot \psi_{V}^{A}\ot W)\co (A\ot \psi_{V}^{A}\ot
\psi_{W}^{A})  $
\item[ ]$\hspace{0.38cm}\co (\nu_{V}\ot \nu_{ W}\ot A) $

\item[ ]$=\nabla_{A\ot V\ot W}\co (\beta_{\nu_{V}}\ot W)\co \beta_{\nu_{W}}  $

\item[ ]$= \nabla_{A\ot V\ot W}\co (\mu_{A}\ot V\ot W)\co (A \ot
(( \beta_{\nu_{V}}\ot W)\co \nu_{W}))$

\item[ ]$=\beta_{\nu_{V\ot W}},$

\end{itemize}

where the first equality follows by (\ref{fi-nab}) for ${\Bbb
A}_{V\ot W}$ and (\ref{falso-idemp2}), the second one follows by the
associativity of $\mu_{A}$, (\ref{wmeas-wcp}) for ${\Bbb A}_{V}$ and (\ref{pre3-wcp}) for $\beta_{\nu_{V}}$ and  $\beta_{\nu_{W}}$, the third one follows by
the left $A$-linearity of $\beta_{\nu_{V}}$ and the last one follows
from the left $A$-linearity of $\nabla_{A\ot V\ot W}$.

The proof for the equality (\ref{pre1-wcp}) is the following:

\begin{itemize}

\item[ ]$\hspace{0.38cm} (\mu_{A}\ot V\ot W)\co (A\ot \sigma_{V\ot
W}^{A})\co (\psi_{V\ot W}^{A}\ot V\ot W)\co (V\ot W\ot \nu_{V\ot W})
$

\item[ ]$= (\mu_{A}\ot V\ot W)\co (A\ot \sigma_{V\ot
W}^{A})\co (\psi_{V\ot W}^{A}\ot V\ot W)\co (V\ot W\ot ((\mu_{A}\ot
V\ot W)\co (A\ot \psi_{V}^{A}\ot W)\co (\nu_{V}\ot \nu_{W}))) $

\item[ ]$=(\mu_{A}\ot V\ot W)\co (\mu_{A}\ot A\ot V\ot W)\co (A\ot \mu_{A}\ot \sigma_{V}^{A}\ot W)\co
(A\ot A\ot \psi_{V}^{A}\ot V\ot W)\co (A\ot \psi_{V}^{A}\ot
\psi_{V}^{A}\ot W) $
\item[ ]$\hspace{0.38cm} \co (A\ot V\ot A\ot V\ot \sigma_{W}^{A})\co
(A\ot V\ot A\ot \tau_{W}^{V}\ot W)\co (\psi_{V}^{A}\ot
\psi_{W}^{A}\ot V\ot W)\co (V\ot \psi_{W}^{A}\ot \psi_{V}^{A} \ot
W)  $
\item[ ]$\hspace{0.38cm} \co (\Delta_{V\ot W}\ot \nu_{V}\ot \nu_{W})  $

\item[ ]$=(\mu_{A}\ot V\ot W)\co (A\ot \sigma_{V}^{A}\ot W)\co
((\psi_{V}^{A}\co (V\ot (\mu_{A}\co (A\ot \mu_{A}))))\ot V\ot W) \co
(V\ot A\ot A\ot \psi_{V}^{A} \ot W) $
\item[ ]$\hspace{0.38cm} \co (V\ot A\ot A\ot V\ot
\sigma_{W}^{A})\co (V\ot A\ot ((A\ot \tau_{W}^{V})\co
(\psi_{W}^{A}\ot V)\co (W\ot \psi_{V}^{A}))\ot W)\co (V\ot
\psi_{W}^{A}\ot  V\ot A\ot W)  $
\item[ ]$\hspace{0.38cm} \co (\Delta_{V\ot W}\ot \nu_{V}\ot \nu_{W}) $

\item[ ]$=(\mu_{A}\ot V\ot W)\co (A\ot \sigma_{V}^{A}\ot W)\co
((\psi_{V}^{A}\co (V\ot (\mu_{A}\co (A\ot \mu_{A}))))\ot V\ot W) \co
(V\ot A\ot A\ot \psi_{V}^{A} \ot W) $
\item[ ]$\hspace{0.38cm}\co (V\ot A\ot A\ot V\ot
\sigma_{W}^{A})\co (V\ot A\ot (((\psi_{V}^{A}\ot W)\co (V\ot
\psi_{W}^{A})\co (\tau_{W}^{V}\ot A)  )\ot W)\co (V\ot
\psi_{W}^{A}\ot  V\ot A\ot W) $
\item[ ]$\hspace{0.38cm} \co (\Delta_{V\ot W}\ot \nu_{V}\ot \nu_{W}) $

\item[ ]$=(\mu_{A}\ot V\ot W)\co (A\ot \sigma_{V}^{A}\ot W)\co
((\psi_{V}^{A}\co (V\ot \mu_{A}))\ot V\ot W)\co (V\ot A\ot
\psi_{V}^{A}\ot W) $
\item[ ]$\hspace{0.38cm}\co (V\ot A\ot V\ot ((\mu_{A}\ot W)\co (A\ot \sigma_{W}^{A})\co
(\psi_{W}^{A}\ot W)\co (W\ot \nu_{W})))\co (V\ot A\ot
\tau_{W}^{V})\co (V\ot \psi_{W}^{A}\ot V)  $
\item[ ]$\hspace{0.38cm}\co (\Delta_{V\ot W}\ot \nu_{V}) $

\item[ ]$=(\mu_{A}\ot V\ot W)\co (A\ot \sigma_{V}^{A}\ot W)\co
(\psi_{V}^{A}\ot V\ot W)\co (V\ot \mu_{A}\ot V\ot W)\co (V\ot A\ot
\psi_{V}^{A}\ot W)\co (V\ot  A\ot V\ot \psi_{W}^{A}) $
\item[ ]$\hspace{0.38cm} \co (V\ot A\ot \tau_{W}^{V}\ot A)\co (V\ot \psi_{W}^{A}\ot V\ot A)\co
(\Delta_{V\ot W}\ot \nu_{V}\ot \eta_{A}) $

\item[ ]$=(\mu_{A}\ot V\ot W)\co (A\ot ((\mu_A\ot V)\co (A\ot \sigma_{V}^{A})\co (\psi_{V}^{A}\ot V)\co
(V\ot \psi_{V}^{A}))\ot W)\co (A\ot V\ot V\ot \psi_{W}^{A})\co (
\psi_{V}^{A}\ot \tau_{W}^{V}\ot A)  $
\item[ ]$\hspace{0.38cm}\co(V\ot  \psi_{W}^{A}\ot V\ot A)\co
(\Delta_{V\ot W}\ot \nu_{V}\ot \eta_{A})   $

\item[ ]$=(\mu_{A}\ot V\ot W)\co (A\ot \psi_{V}^{A}\ot W)\co
(A\ot V\ot (\psi_{W}^{A}\co (W\ot \eta_{A})))\co (\mu_{A}\ot
V\ot W)\co (A\ot \sigma_{V}^{A}\ot W)  $
\item[ ]$\hspace{0.38cm} \co (\psi_{V}^{A}\ot \tau_{W}^{V})\co (V\ot
\psi_{W}^{A}\ot V)\co (\Delta_{V\ot W}\ot \nu_{V})
   $

\item[ ]$= (\mu_{A}\ot V\ot W)\co (A\ot \psi_{V}^{A}\ot W)\co
(A\ot V\ot (\psi_{W}^{A}\co (W\ot \eta_{A})))\co
(\psi_{V}^{A}\ot W)\co (V\ot \psi_{W}^{A})\co (\Delta_{V\ot W}\ot
 \eta_{A})
 $

\item[ ]$=\nabla_{A\ot V\ot W}\co (\eta_{A}\ot V\ot W).
 $

\end{itemize}

In this proof, the first equality follows by (\ref{c11}) for ${\Bbb
A}_{V\ot W}$, the second one follows by (\ref{wmeas-wcp})  for ${\Bbb
A}_{V}$, ${\Bbb A}_{W}$, the associativity of $\mu_{A}$ and the
twisted condition for ${\Bbb A}_{V}$, the third one follows by
(\ref{wmeas-wcp})  for ${\Bbb A}_{V}$ and the associativity of
$\mu_{A}$. In the fourth one we applied (i) of the definition of
twisting morphism and the fifth one is a consequence of
(\ref{wmeas-wcp})  for ${\Bbb A}_{V}$. The sixth one follows by
(\ref{pre1-wcp}) for ${\Bbb A}_{W}$, the seventh one follows by
(\ref{wmeas-wcp})  for ${\Bbb A}_{V}$ and the ssociativity of
$\mu_{A}$, the eight one relies on the associativity of $\mu_{A}$ and the
twisted condition for ${\Bbb A}_{V}$. Finally, the ninth one is a
consequence of (\ref{pre-1}) and the last one follows by  (\ref{wmeas-wcp})  for
${\Bbb A}_{V}$, ${\Bbb A}_{W}$.

On the other hand, the proof for the identity (\ref{pre2-wcp}) is

\begin{itemize}

\item[ ]$\hspace{0.38cm} (\mu_{A}\ot V\ot W)\co (A\ot \sigma_{V\ot
W}^{A})\co  (\nu_{V\ot W}\ot V\ot W) $

\item[ ]$= (\mu_{A}\ot V\ot W)\co (A\ot \sigma_{V\ot
W}^{A})\co (((\mu_{A}\ot V\ot W)\co (A\ot \psi_{V}^{A}\ot W)\co
(\nu_{V}\ot \nu_{W}))\ot V\ot W) $

\item[ ]$=(\mu_{A}\ot V\ot W)\co  ((\mu_{A}\co (A\ot \mu_{A}))\ot \sigma_{V}^{A}\ot W)\co
(A\ot A\ot \psi_{V}^{A}\ot V\ot W)\co (A\ot A\ot V\ot
\psi_{V}^{A}\ot W)  $
\item[ ]$\hspace{0.38cm} \co (A\ot A\ot V\ot V\ot \sigma_{W}^{A})\co
(A\ot \psi_{V}^{A}\ot \tau_{W}^{V}\ot W)\co (\nu_{V}\ot \nu_{W}\ot
V\ot W)    $

\item[ ]$=(\mu_{A}\ot V\ot W)\co  (\mu_{A}\ot \sigma_{V}^{A}\ot W)\co
(A\ot \psi_{V}^{A}\ot V\ot W)\co (\nu_{V}\ot ((\mu_{A}\ot V\ot W)\co
(A\ot \psi_{V}^{A}\ot W)\co (A\ot V\ot \sigma_{W}^{A})  $
\item[ ]$\hspace{0.38cm} \co (A\ot \tau_{W}^{V}\ot W)\co (\nu_{W}\ot V\ot W) $

\item[ ]$=(\mu_{A}\ot V\ot W)\co  (\mu_{A}\ot \sigma_{V}^{A}\ot W)\co
(A\ot \psi_{V}^{A}\ot V\ot W)\co (\nu_{V}\ot (\nabla_{A\ot V\ot
W}\co (\eta_{A}\ot V\ot W))) $

\item[ ]$=(\mu_{A}\ot V\ot W)\co  (A\ot \psi_{V}^{A}\ot W)\co (A\ot V\ot
(\psi_{W}^{A}\co (W\ot \eta_{A})))\co (((\mu_{A}\ot V)\co (A\ot
\sigma_{V}^{A})\co (\nu_{V}\ot V))\ot W)   $
\item[ ]$\hspace{0.38cm}\co \Delta_{V\ot W}$

\item[ ]$=(\mu_{A}\ot V\ot W)\co  (A\ot \psi_{V}^{A}\ot W)\co (A\ot V\ot
(\psi_{W}^{A}\co (W\ot \eta_{A})))\co ((\nabla_{A\ot V}\co
(\eta_{A}\ot V))\ot W)\co \Delta_{V\ot W}  $

\item[ ]$=\nabla_{A\ot V\ot W}\co (\eta_{A}\ot V\ot W).
 $

\end{itemize}

The first equality follows by (\ref{aw1}) for ${\Bbb A}_{V\ot W}$,
the second one follows by the associativity of $\mu_{A}$ and the twisted
condition for ${\Bbb A}_{V}$, the third one follows by (\ref{wmeas-wcp})
for ${\Bbb A}_{V}$, the fourth one relies on (\ref{pre-2}), the fifth one
follows by the twisted condition for ${\Bbb A}_{V}$ and the associativity of $\mu_{A}$, the sixth one is a consequence of (\ref{pre2-wcp}) for ${\Bbb A}_{V}$ and finally, the last one follows by
(\ref{wmeas-wcp}) for ${\Bbb A}_{V}$.

\end{dem}

\begin{nota}
{\rm Note that we can obtain similar results about the iteration process if we work with quadruples $\;_{V}{\Bbb A}=(V,A,\psi_{A}^{V}, \sigma_{A}^{V})$ where $\psi_{A}^{V}:A\ot V\rightarrow V\ot A$ and $\sigma_{A}^{V}: V\ot V\rightarrow A\ot V$ satisfy the suitable conditions that define a weak crossed product on $V\ot A$.
}
\end{nota}

\section{Some examples}

The aim of this section is to provide some examples of the iteration process introduced in the previous one.

\begin{exemplo}
\label{ej1} {\rm Let ${\mathcal C}$ be a category. The category of endofunctors of
${\mathcal C}$ is a strict monoidal category with the composition of
functors, denoted by $\circledcirc$,  as the product and the
identity functor as the unit. We denote this category by
$End({\mathcal C})$. The morphisms in $End({\mathcal C})$ are
natural transformations between endofunctors and we denote the
composition (the vertical composition) of these morphisms by
$\circ$. The tensor product of morphisms in $End({\mathcal C})$ is
defined by the horizontal composition of natural transformations and
in this paper is denoted by the same symbol used for the composition
of functors (see \cite{Macl} for the details of the horizontal and
vertical compositions). Given objects $T$, $S$, $H$ and a morphism
$\tau:S\rightarrow H$, we write $T\circledcirc \tau$ for
$id_{T}\circledcirc \tau$ and $\tau\circledcirc T$ for
$\tau\circledcirc id_{T}$ where $id_{T}$ denotes the identity
morphism for the object $T$.

A monad on ${\mathcal C}$ consists of a endofunctor $S:{\mathcal
C}\rightarrow {\mathcal C}$ together with two natural
transformations $\eta_{S}:id_{{\mathcal C}}\rightarrow S$ (where
$id_{{\mathcal C}}$ denotes the identity functor on ${\mathcal C}$)
and $\mu_{S}:S^{2}=S\circledcirc S\rightarrow S$. These are required
to fulfill the following conditions
\begin{equation}
\label{alg1}
 \mu_{S}\circ (S\circledcirc \eta_{S})=\mu_{S}\circ
(\eta_{S}\circledcirc S)=id_{S},
\end{equation}
\begin{equation}
\label{alg2}\mu_{S}\circ (S\circledcirc \mu_{S})=\mu_{S}\circ
(\mu_{S}\circledcirc S).
\end{equation}

Then, a monad on ${\mathcal C}$ can alternatively be defined as a
monoid in the  strict monoidal category $End({\mathcal
C})$.

The notion of wreath was introduced by Lack and Street in \cite{LS}.
A monad $S$ in ${\mathcal C}$ is a wreath if there exist an object
in $T\in End({\mathcal C})$ and morphisms in $End({\mathcal C})$,
$\lambda:T\circledcirc S\rightarrow S\circledcirc T$,
$\tau:id_{{\mathcal C}}\rightarrow S \circledcirc T$ and $v:T
\circledcirc T\rightarrow S \circledcirc T$ satisfying the following
conditions:
\begin{equation}
\label{W1}(\mu_{S}\circledcirc T)\circ (S\circledcirc \lambda)\circ
(\lambda\circledcirc S)=\lambda\circ (T\circledcirc \mu_{S}),
\end{equation}
\begin{equation}
\label{W2}\lambda\circ (T\circledcirc \eta_{S})=\eta_{S}\circledcirc
T,\end{equation}
\begin{equation}
\label{W3}(\mu_{S}\circledcirc T)\circ (S \circledcirc
\tau)=(\mu_{S}\circledcirc T)\circ (S\circledcirc \lambda)\circ
(\tau\circledcirc S),\end{equation}
\begin{equation}
\label{W4}(\mu_{S}\circledcirc T)\circ (S\circledcirc v)\circ (\lambda
\circledcirc T)\circ (T \circledcirc \lambda)=(\mu_{S}\circledcirc
T)\circ (S \circledcirc \lambda)\circ (v\circledcirc S),
\end{equation}
\begin{equation}
\label{W5}(\mu_{S}\circledcirc T)\circ (S\circledcirc v)\circ
(v\circledcirc T)=(\mu_{S}\circledcirc T)\circ (S\circledcirc
v)\circ (\lambda\circledcirc T)\circ (T\circledcirc
v),\end{equation}
\begin{equation}
\label{W6}(\mu_{S}\circledcirc T)\circ (S\circledcirc v)\circ
(\tau\circledcirc T)=\eta_{S}\circledcirc T=(\mu_{S}\circledcirc
T)\circ (S\circledcirc v)\circ (\lambda\circledcirc T)\circ
(T \circledcirc \tau).
\end{equation}

If we put $\psi_{T}^{S}=\lambda$ and $\sigma_{T}^{S}=v$, we obtain
that ${\Bbb S}_{T}=(S,T, \psi_{T}^{S}, \sigma_{T}^{S})$ is a quadruple satisfying
(\ref{wmeas-wcp}), (\ref{twis-wcp}) and (\ref{cocy2-wcp}) where the
associated idempotent defined in (\ref{idem-wcp}) is
$\nabla_{S\circledcirc T}=id_{S\circledcirc T}$ because $\lambda$
satisfies the identity (\ref{W2}). Then, the product induced by a
wreath (wreath product) defined by
$$\mu_{S\circledcirc T}=(\mu_{S}\circledcirc T)\circ
(\mu_{S}\circledcirc v)\circ (S\circledcirc \lambda\circledcirc T)$$
is the one defined in (\ref{prod-todo-wcp}) and it is associative
because satisfies (iv) (twisted condition) and (v) (cocycle
condition). Then $S\circledcirc T$ is a monad with  unit
$\eta_{S\circledcirc T}=\tau$.

Note that, in this case we do not need that every  idempotent splits
because the associated idempotent $\nabla_{S\circledcirc
T}=id_{S\circledcirc T}$. Therefore wreath products are examples of
weak crossed products with trivial idempotent.

An example of wreath products cames from the notion of distributive
law introduced by Beck in \cite{Beck} (see also \cite{Street-FTM}).
Suppose that $T$ and $S$ are two monads on ${\mathcal C}$. A
distributive law of the monad $S$ over the monad $T$ is a natural
transformation
$$\lambda:T\circledcirc S\rightarrow S\circledcirc T$$
such that
\begin{equation}
\label{DL1}
 \lambda\circ (\mu_{T}\circledcirc S)=
(S\circledcirc \mu_{T})\circ (\lambda\circledcirc T)\circ
(T\circledcirc \lambda),
\end{equation}
\begin{equation}
\label{DL2}\lambda\circ (\eta_{T}\circledcirc S)=S\circledcirc
\eta_{T},
\end{equation}
\begin{equation}
\label{DL3} \lambda\circ (T\circledcirc\mu_{S})=
(\mu_{S}\circledcirc T)\circ (S\circledcirc \lambda)\circ (\lambda
\circledcirc S),
\end{equation}
\begin{equation}
\label{DL4} \lambda\circ (T\circledcirc
\eta_{S})=\eta_{S}\circledcirc T.
\end{equation}

Then, if  $\tau=\eta_{S}\circledcirc\eta_{T}$ and
$v=\eta_{S}\circledcirc\mu_{T}$ we obtain a wreath for the monad $S$
and also a weak crossed product associated to the quadruple ${\Bbb S}_{T}=(S,T, \psi_{T}^{S}, \sigma_{T}^{S})$ where $\psi_{T}^{S}=\lambda$,
$\sigma_{T}^{S}=v$ and
$$\mu_{S\circledcirc T}=(\mu_{S}\circledcirc \mu_{T})\circ (S\circledcirc
\lambda \circledcirc T).$$

Suppose that $S$, $T$ and $D$ are monads in ${\mathcal C}$ such that
there exists the following distributive laws between them
$$\lambda_{1}:T\circledcirc S\rightarrow S\circledcirc T,\;\;\;\;
\lambda_{2}:D\circledcirc T\rightarrow T\circledcirc D, \;\;\;\;
\lambda_{3}:D\circledcirc S\rightarrow S\circledcirc D,$$ satisfying
the  compatibility identity (called the Yang-Baxter relation or the
hexagon equation)
\begin{equation}
\label{YB-Comp} (S\circledcirc \lambda_{2})\circ
(\lambda_{3}\circledcirc T)\circ (D\circledcirc
\lambda_{1})=(\lambda_{1}\circledcirc D)\circ (T\circledcirc
\lambda_{3})\circ (\lambda_{2}\circledcirc S).
\end{equation}

Then, under these conditions we have two quadruples
$${\Bbb S}_{T}=(S, T, \psi_{T}^{S}=\lambda_{1},
\sigma_{T}^{S}=\eta_{S}\circledcirc\mu_{T}),$$
$${\Bbb S}_{D}=(S, D, \psi_{D}^{S}=\lambda_{3},
\sigma_{D}^{S}=\eta_{S}\circledcirc\mu_{D}),$$ satisfying
(\ref{wmeas-wcp}), (\ref{twis-wcp}), (\ref{cocy2-wcp}). If we put
$\Delta_{T\circledcirc D}=id_{T\circledcirc D}$ as a link morphism
(note that in this case the equalities (\ref{sigma1}),
(\ref{sigma2}) and (\ref{sigma3}) are trivial) and
$\tau_{D}^{T}=\lambda_{2}$ we have that the condition (i) of
Definition \ref{wcp-def} holds because we assume (\ref{YB-Comp}). On
the other hand, the condition (ii) of the same Definition also holds
because:

\begin{itemize}

\item[ ]$\hspace{0.38cm}(\mu_{S}\circledcirc T\circledcirc D)\co (S\circledcirc
\sigma_{T}^{S}\circledcirc D)\co (\psi_{T}^{S}\circledcirc
\tau_{D}^{T})\circledcirc (T\circledcirc \sigma_{D}^{S}\circledcirc
T)\co (\tau_{D}^{T}\circledcirc D\circledcirc T)$

\item [ ]$=(\eta_{S}\circledcirc ((\mu_{T}\circledcirc\mu_{D})\circ
(T\circledcirc \lambda_{2}\circledcirc D)\circ
(\lambda_{2}\circledcirc\lambda_{2})))$

\item[ ]$=(\mu_{S}\circledcirc T\circledcirc D)\co (S\circledcirc \psi_{T}^{S}\circledcirc
D)\co (S\circledcirc T\circledcirc \sigma_{D}^{S})\co (S\circledcirc
\tau_{D}^{T} \circledcirc D)\co (\psi_{D}^{S}\circledcirc
T\circledcirc D)\co (D\circledcirc \sigma_{T}^{S}\circledcirc D)\co
 $
\item[ ]$\hspace{0.38cm}\co (D\circledcirc T\circledcirc \tau_{D}^{T}).$

\end{itemize}

Therefore, $\tau_{D}^{T}=\lambda_{2}$ is a twisting morphism between
the quadruples ${\Bbb S}_{T}$ and ${\Bbb S}_{D}$. As a consequence,
by Lemma \ref{iden-ite} and Theorem \ref{prince}, the quadruple
$${\Bbb S}_{T\circledcirc D}=(S, T\circledcirc D,
\psi_{T\circledcirc D}^{S}, \sigma_{T\circledcirc D}^{S}),$$ where
$$\psi_{T\circledcirc D}^{S}=
(\psi_{T}^{S}\circledcirc D)\co (T\circledcirc
\psi_{D}^{S})=(\lambda_{1}\circledcirc D)\co (T\circledcirc
\lambda_{3})$$ and
$$ \sigma_{T\circledcirc D}^{S}=
(\mu_{S}\circledcirc T\circledcirc D)\co (S\circledcirc
\psi_{T}^{S}\circledcirc D)\co (\sigma_{T}^{S}\circledcirc
\sigma_{D}^{S})\co (T\circledcirc \tau_{D}^{T}\ot D)=
\eta_{S}\circledcirc ( (\mu_{T}\circledcirc \mu_{D})\circ
(T\circledcirc \lambda_{2}\circledcirc D)),$$ satisfies the
equalities (\ref{twis-wcp}) and (\ref{cocy2-wcp}). Then, the pair,
$(S\circledcirc T\circledcirc  D, \mu_{S\circledcirc T\circledcirc
D})$ is the iterated weak crossed of $(S\circledcirc T,
\mu_{S\circledcirc T})$ and $(S\circledcirc  D, \mu_{S\circledcirc
D})$ with associated product
$$\mu_{S\circledcirc T\circledcirc D}=(\mu_S\circledcirc T\circledcirc D)\co (\mu_S\circledcirc
\sigma_{T\circledcirc D}^{S})\co (S\circledcirc \psi_{T\circledcirc
D}^{S}\circledcirc T \circledcirc D)=$$ $$(\mu_{S}\circledcirc
\mu_{T}\circledcirc\mu_{D})\circ (S\circledcirc
((\lambda_{1}\circledcirc\lambda_{2})\circ
(T\circledcirc\lambda_{3}\circledcirc T))\circledcirc D).$$

In this case the preunits are units. The object $S\circledcirc T\circledcirc
D$ is a monad with unit
$$\eta_{S\circledcirc T\circledcirc D}=\eta_{S}\circledcirc \eta_{T}\circledcirc \eta_{D}$$
because $S\circledcirc T$ and $S\circledcirc D$ are also monads with unit
$\eta_{S\circledcirc T}=\eta_{S}\circledcirc\eta_{T}$ and
$\eta_{S\circledcirc D}=\eta_{S}\circledcirc\eta_{D}$ respectively.
Therefore, (\ref{pre-1}) and (\ref{pre-2}) holds and the morphism
$\nu_{T\circledcirc S}$ defined in (\ref{iterated-preunit}) is
$\eta_{S\circledcirc T\circledcirc D}$.

For example, if ${\mathcal C}$ is a strict monoidal category and
$A$, $B$ are monoids in ${\mathcal C}$ the twisted tensor product
of algebras introduced in \cite{CSV}, \cite{TAM} is an example weak
crossed product associated to a wreath for the monad $S=A\otimes -$.
In this case $T=B\otimes -$ and $\lambda=R\otimes -$ where
$R:B\otimes A\rightarrow A\otimes B$  is the twisting morphism.
Furthermore, the natural transformation $\lambda=R\otimes -$ is a
distributive law of the monad $S=A\otimes -$ over the monad
$T=B\otimes -$ if and only if $R:B\otimes A\rightarrow A\otimes B$
is a unital twisting morphism. Suppose that $A$, $B$ and $C$ are
monoids, let
$$ R_1:B\ot A\rightarrow A\ot B,\;\;
R_2:C\ot B\rightarrow B\otimes C,\;\; R_3:C\ot A\rightarrow A\ot C,
$$
unital twisting morphisms, and consider  the monads $S=A\ot - $,
$T=B\ot -$, $D=C\ot -$, the induced quadruples ${\Bbb S}_{T}$,
${\Bbb S}_{D}$ and the twisting morphism $\tau_{D}^{T}=R_{2}\otimes
-$. Then the iterated product defined in Theorem 2.1 of
\cite{Pascual-Javi2} is the  one associated to the  quadruple
$${\Bbb S}_{T\circledcirc D}=(S, T\circledcirc D,
\psi_{T\circledcirc D}^{S}, \sigma_{T\circledcirc D}^{S})$$ when we
apply the functors in the unit object of the category.
}
\end{exemplo}

\begin{exemplo}
\label{ej2} {\rm In this case we assume that ${\mathcal C}$ is a
category where every idempotent morphism splits. As in the previous
example we work in the category $End({\mathcal C})$ and it  is easy
to show that every idempotent morphism splits in  $End({\mathcal
C})$ because every idempotent morphism splits in ${\mathcal C}$.

Given to monads $S$ and $T$,  the notion of weak distributive law of
the monad $S$ over the monad $T$ was introduced by Ross Street in
\cite{Street-WDL} as follows. It consists of a natural
transformation
$$\lambda:T\circledcirc S\rightarrow S\circledcirc T$$
such that satisfies (\ref{DL1}), (\ref{DL3}) and

\begin{equation}
\label{WDL1} \lambda\circ (\eta_{T}\circledcirc
S)=(\mu_{S}\circledcirc T)\circ (S\circledcirc (\lambda\circ
(\eta_{T}\circledcirc \eta_{S}))),
\end{equation}
\begin{equation}
\label{WDL2}\lambda\circ (T\circledcirc \eta_{S})=(S\circledcirc
\mu_{T})\circ ((\lambda\circ (\eta_{T}\circledcirc
\eta_{S}))\circledcirc S).
\end{equation}

In this definition the axioms (\ref{WDL1}) and (\ref{WDL2}) can be
replaced for the identity [\cite{Street-WDL}, Proposition 2.2]:
\begin{equation}
\label{idem=idem} (S\circledcirc \mu_{T})\circ ((\lambda\circ
(\eta_{T}\circledcirc S)\circledcirc T)=(\mu_{S}\circledcirc T)\circ
(S\circledcirc (\lambda \circ (T\circledcirc\eta_{S}))).
\end{equation}

For a weak distributive law, the weak wreath product is defined by
$$\mu_{S\circledcirc T}=(\mu_{S}\circledcirc\mu_{T})\circ
(S\circledcirc\lambda \circledcirc T).$$

The same set of axioms for monoids in category of modules over a
commutative ring can be found in \cite{Caen}. Then, the conditions
used in \cite{Caen} define a weak wreath product associated to
monads induced by monoids.

It follows by (\ref{DL1}) and (\ref{DL3}) that $\mu_{S\circledcirc
T}$ is an associative product but possibly without unity. In any
case, if we take the quadruple
$${\Bbb S}_{T}=(S, T, \psi_{T}^{S}=\lambda,
\sigma_{T}^{S}=(S\circledcirc \mu_{T})\circ ((\lambda\circ
(T\circledcirc \eta_{S}))\circledcirc T)),$$ we obtain that ${\Bbb
S}_{T}$ satisfies (\ref{wmeas-wcp}), (\ref{twis-wcp}),
(\ref{cocy2-wcp}) and (\ref{idemp-sigma-inv}). The associated idempotent defined in
(\ref{idem-wcp}) is
$$\nabla_{S\circledcirc T}=(\mu_{S}\circledcirc
T)\circ (S\circledcirc (\lambda \circ (T\circledcirc\eta_{S}))).$$

Then, the weak wreath product defined by the weak distributive law
is the one induced by the quadruple ${\Bbb S}_{T}$. Therefore, every
weak wreath product is a weak crossed product. In this setting the
morphism $\nu_{T} =\nabla_{S\circledcirc T}\circ (\eta_{S}\circledcirc
\eta_{T})$ is a preunit and $S\times T $ is a monoid with  unit
$\eta_{S\times T}=p_{S\circledcirc T}\circ\nu_{T}$ (see also
[\cite{mra-preunit}, Example 3.16]).

Note that the equality (\ref{idem=idem}) implies that
\begin{equation}
\label{equ-idem} \sigma_{T}^{S}=(S\circledcirc \mu_{T})\circ
((\nabla_{S\circledcirc T}\circ (\eta_{S}\circledcirc T)\circledcirc
T)=\nabla_{S\circledcirc T}\circ
(\eta_{S}\circledcirc\mu_{T})=\lambda\circ (\mu_{T}\circledcirc
\eta_{S}).
\end{equation}

Suppose that $S$, $T$ and $D$ are monads in ${\mathcal C}$ such that
there exists  tree weak distributive laws between them
$$\lambda_{1}:T\circledcirc S\rightarrow S\circledcirc T,\;\;\;\;
\lambda_{2}:D\circledcirc T\rightarrow T\circledcirc D, \;\;\;\;
\lambda_{3}:D\circledcirc S\rightarrow S\circledcirc D,$$ satisfying
 the Yang-Baxter relation (\ref{YB-Comp}). Then, under these conditions we have two quadruples
$${\Bbb S}_{T}=(S, T, \psi_{T}^{S}=\lambda_{1},
\sigma_{T}^{S}=(S\circledcirc \mu_{T})\circ ((\lambda_{1}\circ
(T\circledcirc \eta_{S}))\circledcirc T)),$$
$${\Bbb S}_{D}=(S, D, \psi_{D}^{S}=\lambda_{3},
\sigma_{D}^{S}=(S\circledcirc \mu_{D})\circ ((\lambda_{3}\circ
(D\circledcirc \eta_{S}))\circledcirc D)),$$ satisfying
(\ref{wmeas-wcp}), (\ref{twis-wcp}), (\ref{cocy2-wcp}) and (\ref{idemp-sigma-inv}).  If we put $\Delta_{T\circledcirc
D}=\nabla_{T\circledcirc D}$ we obtain a link morphism. Indeed,  we have that
(\ref{falso-idemp}) holds because

\begin{itemize}

\item[ ]$\hspace{0.38cm} (S\circledcirc \nabla_{T\circledcirc D})\co \psi_{T\circledcirc D}^{S}$

\item [ ]$= (S\circledcirc \mu_{T}\circledcirc\mu_{D})\co (\lambda_{1}\circledcirc\lambda_{2}\circledcirc D)\co (T\circledcirc \lambda_{3}\circledcirc \lambda_{2})\co ((\lambda_{2}\co (\eta_{D}\circledcirc T))\circledcirc \lambda_{3}\circledcirc\eta_{T}) $

\item[ ]$= (S \circledcirc ((\mu_{T}\circledcirc\mu_{D})\co (T\circledcirc\lambda_{2}\circledcirc D)\co (\lambda_{2}\circledcirc T\circledcirc D)))\co (\lambda_{3}\circledcirc T\circledcirc T\circledcirc D)\co (D\circledcirc \lambda_{1}\circledcirc \lambda_{2})\co (\eta_{D}\circledcirc D\circledcirc \lambda_{3}\circledcirc \eta_{T})$

\item [ ]$= (S\circledcirc T\circledcirc \mu_{D})\co (S\circledcirc\lambda_{2} \circledcirc D)\co (\lambda_{3}\circledcirc \nabla_{T\circledcirc D})\co (D\circledcirc \lambda_{1}\circledcirc D)\co (\eta_{D}\circledcirc T\circledcirc \lambda_{3}) $

\item [ ]$= (S\circledcirc T\circledcirc \mu_{D})\co (((S\circledcirc\lambda_{2})\co (\lambda_{3}\circledcirc T)\co (D\circledcirc\lambda_{1}))\circledcirc D)\co (\eta_{D}\circledcirc T\circledcirc \lambda_{3}) $

\item [ ]$= (\lambda_{1} \circledcirc \mu_{D})\co (T\circledcirc \lambda_{3}\circledcirc D)\co ((\lambda_{2}\co (\eta_{D}\circledcirc T))\circledcirc \lambda_{3})$

\item[ ]$= \psi_{T\circledcirc D}^{S}  $

\end{itemize}

In the last equalities, the first one follows by (\ref{DL1}) for $\lambda_{3}$ and $\lambda_{2}$, the second one follows by (\ref{YB-Comp}), the third one follows by (\ref{DL3}) for $\lambda_{2}$, the fourth one follows by
\begin{equation}
\label{new-nabla}
(T\circledcirc \mu_{D})\co (\lambda_{2}\circledcirc D)\co (D\circledcirc \nabla_{T\circledcirc D})=(T\circledcirc \mu_{D})\co (\lambda_{2}\circledcirc D),
\end{equation}
the fifth one relies on (\ref{YB-Comp}) and the last one is a consequence of (\ref{DL1}) for $\lambda_{3}$.

The equality (\ref{falso-idemp2}) follows by

\begin{itemize}

\item[ ]$\hspace{0.38cm} \nabla_{S\circledcirc T\circledcirc D}\co (\lambda_{1}\circledcirc D)\co (T\circledcirc \lambda_{3})$

\item [ ]$= (((\mu_{S} \circledcirc \mu_{T})\co (S\circledcirc\lambda_{1}\circledcirc T))\circledcirc D)\co (\lambda_{1}\circledcirc ((\lambda_{1}\circledcirc D)\co 
    (T\circledcirc \lambda_{3})\co (\lambda_{2}\circledcirc S)))\co (T\circledcirc \lambda_{3}\circledcirc\eta_{T}\circledcirc \eta_{S})$

\item[ ]$= (((\mu_{S} \circledcirc \mu_{T})\co (S\circledcirc\lambda_{1}\circledcirc T))\circledcirc D)\co (\lambda_{1}\circledcirc ((S\circledcirc \lambda_{2})\co (\lambda_{3} \circledcirc T)\co (D\circledcirc \lambda_{1})))\co (T\circledcirc \lambda_{3}\circledcirc\eta_{T}\circledcirc \eta_{S})$

\item [ ]$= (S\circledcirc\mu_{T}\circledcirc D)\co (\lambda_{1}\circledcirc \lambda_{2})\co (T\circledcirc\lambda_{3}\circledcirc T)\co (T\circledcirc D \circledcirc (\nabla_{S\circledcirc T}\co (S\circledcirc \eta_{T}))) $

\item [ ]$= (S\circledcirc\mu_{T}\circledcirc D)\co (\lambda_{1}\circledcirc T\circledcirc D)\co (T\circledcirc (((S\circledcirc \lambda_{2})\co (\lambda_{3} \circledcirc T)\co (D\circledcirc \lambda_{1}))\co (D\circledcirc \eta_{T}\circledcirc S))) $

\item [ ]$= (S\circledcirc\mu_{T}\circledcirc D)\co (\lambda_{1}\circledcirc T\circledcirc D)\co (T\circledcirc (((\lambda_{1}\circledcirc D)\co
    (T\circledcirc \lambda_{3})\co (\lambda_{2}\circledcirc S))\co (D\circledcirc \eta_{T}\circledcirc S)))$

\item[ ]$=\psi_{T\circledcirc D}^{S}    $

\end{itemize}

where, the first and the sixth equalities follow by (\ref{DL1}) for $\lambda_{1}$, the second and the fifth ones follow by (\ref{YB-Comp}), the third relies on (\ref{DL3}) for $\lambda_{1}$ and $\lambda_{3}$ and, finally, the fourth one follows by (\ref{idem=idem}).

On the other hand, if 
$\tau_{D}^{T}=\lambda_{2}$ we obtain that the condition (i) of
Definition \ref{wcp-def} holds by (\ref{YB-Comp}). Moreover,  condition (ii) of the same Definition also holds because we have the  following:

\begin{itemize}

\item[ ]$\hspace{0.38cm}(\mu_{S}\circledcirc T\circledcirc D)\co (S\circledcirc
\sigma_{T}^{S}\circledcirc D)\co (\psi_{T}^{S}\circledcirc
\tau_{D}^{T})\circledcirc (T\circledcirc \sigma_{D}^{S}\circledcirc
T)\co (\tau_{D}^{T}\circledcirc D\circledcirc T)$

\item [ ]$= (S \circledcirc\mu_{T}\circledcirc D)\co
(\lambda_{1}\circledcirc (\lambda_{2}\co (\mu_{D}\circledcirc
T)))\co (((T\circledcirc \lambda_{3})\co
(\lambda_{2}\circledcirc\eta_{S}))\circledcirc D\circledcirc T) $

\item[ ]$=(S \circledcirc  ((\mu_{T}\circledcirc \mu_{D})\co
(T\circledcirc\lambda_{2} \circledcirc D)\co
(\lambda_{2}\circledcirc  T\circledcirc D)))\co
(((\lambda_{3}\circledcirc  T)\co (D \circledcirc (\lambda_{1}\co
(T\circledcirc \eta_{S}))))\circledcirc \lambda_{2})$

\item[ ]$=(S\circledcirc T \circledcirc \mu_{D})\co (((S\circledcirc
\lambda_{2})\co (\lambda_{3}\circledcirc \mu_{T})\co (D\circledcirc
(\lambda_{1}\co (T\circledcirc \eta_{S}))\circledcirc
T))\circledcirc D)\co (D\circledcirc T \circledcirc \lambda_{2})$

\item[ ]$=(S\circledcirc T \circledcirc \mu_{D})\co (((S\circledcirc
\lambda_{2})\co (\lambda_{3}\circledcirc T)\co (D\circledcirc
(\lambda_{1}\co (\mu_{T}\circledcirc \eta_{S}))))\circledcirc
 D)\co (D\circledcirc T \circledcirc \lambda_{2})$

\item[ ]$=(\lambda_{1}\circledcirc \mu_{D})\co (T \circledcirc
\lambda_{3}\circledcirc D)\co ((\lambda_{2}\co (D\circledcirc
\mu_{T}))\circledcirc \eta_{S}\circledcirc D)\co (D\circledcirc
T\circledcirc \lambda_{2})$

\item[ ]$=(((S\circledcirc  \mu_{T})\co (\lambda_{1} \circledcirc
T)) \circledcirc  D)\co (T\circledcirc  ((\lambda_{1}\circledcirc
\mu_{D})\co (T\circledcirc \lambda_{3}\circledcirc D)\co
(\lambda_{2}\circledcirc \eta_{S}\circledcirc D)))\co
(\lambda_{2}\circledcirc \lambda_{2})$

\item[ ]$=(S\circledcirc \mu_{T}\circledcirc \mu_{D})\co
(((\lambda_{1}\circledcirc \lambda_{2})\co
(T\circledcirc\lambda_{3}\circledcirc T)\co (\lambda_{2}\circledcirc
(\lambda_{1}\co (T\circledcirc \eta_{S}))))\circledcirc D)\co
(D\circledcirc T\circledcirc \lambda_{2})$

\item[ ]$=(S\circledcirc ((\mu_{T}\circledcirc\mu_{D})\co
(T\circledcirc \lambda_{2}\circledcirc D)\co
(\lambda_{2}\circledcirc T\circledcirc D)))\co
(\lambda_{3}\circledcirc T\circledcirc T \circledcirc D)\co
(D\circledcirc (\lambda_{1}\co (T\circledcirc \mu_{S}))\circledcirc
T\circledcirc D) $
\item[ ]$\hspace{0.38cm}\co (D\circledcirc T\circledcirc
\eta_{S}\circledcirc (((\lambda_{1}\circ (T\circledcirc
\eta_{S}))\circledcirc D)\co \lambda_{2}))$

\item[ ]$=(\mu_{S}\circledcirc  ((\mu_{T}\circledcirc  \mu_{D})\co
(T\circledcirc  \lambda_{2} \circledcirc D)\co
(\lambda_{2}\circledcirc  T\circledcirc  D)))\co (S\circledcirc
\lambda_{3}\circledcirc  T\circledcirc T\circledcirc  D)\co
(\lambda_{3}\circledcirc \lambda_{1}\circledcirc T\circledcirc
D)$
\item[ ]$\hspace{0.38cm}\co (D\circledcirc (\lambda_{1}\co
(T\circledcirc \eta_{S}))\circledcirc (((\lambda_{1}\co
(T\circledcirc \eta_{S}))\circledcirc D)\co \lambda_{2}))$

\item[ ]$=(\mu_{S}\circledcirc T\circledcirc \mu_{D})\co
(S\circledcirc S\circledcirc \lambda_{2}\circledcirc D)\co
(S\circledcirc \lambda_{3}\circledcirc T\circledcirc D)\co 
(\lambda_{3}\circledcirc (\lambda_{1}\co (\mu_{T}\circledcirc
\eta_{S})))\circledcirc D)\co (D\circledcirc (\lambda_{1}\co
(T\circledcirc \eta_{S}))\circledcirc \lambda_{2})$

\item[ ]$=(\mu_{S}\circledcirc T\circledcirc D)\co (S\circledcirc \psi_{T}^{S}\circledcirc
D)\co (S\circledcirc T\circledcirc \sigma_{D}^{S})\co (S\circledcirc
\tau_{D}^{T} \circledcirc D)\co (\psi_{D}^{S}\circledcirc
T\circledcirc D)\co (D\circledcirc \sigma_{T}^{S}\circledcirc D)\co
 $
\item[ ]$\hspace{0.38cm}\co (D\circledcirc T\circledcirc \tau_{D}^{T}).$

\end{itemize}

In the last equalities, the fist one follows by (\ref{DL3}) for
$\lambda_{1}$ and 
$$(S\circledcirc \mu_{T})\co (\lambda_{1}\circledcirc T)\co (T\circledcirc \nabla_{S\circledcirc T})=(S\circledcirc \mu_{T})\co (\lambda_{1}\circledcirc T),$$ the second one follows by
(\ref{DL1}) for $\lambda_{2}$ and (\ref{YB-Comp}), the third one follows by
(\ref{DL3}) for $\lambda_{2}$ and the fourth one is a consequence of 
(\ref{equ-idem}) for $\lambda_{1}$. In the fifth one we used (\ref{YB-Comp}) and the
sixth one relies on (\ref{DL3}) for $\lambda_{2}$ and
(\ref{DL1}) for $\lambda_1$. The seventh one follows by
(\ref{YB-Comp}), the eighth one follows by $\lambda_{1}\co (T\circledcirc \eta_{S})=\nabla_{S\circledcirc T}\co (\eta_{S}\circledcirc T)$ and
(\ref{YB-Comp}). Finally, in the ninth one we used (\ref{DL3}) for
$\lambda_{1}$ and $\lambda_{3}$, the tenth one follows by
(\ref{DL3}) for $\lambda_{2}$ and (\ref{DL1}) for $\lambda_{1}$ and
the last one follows by (\ref{YB-Comp}).

Therefore, $\tau_{D}^{T}=\lambda_{2}$ is a twisting morphism between
the quadruples ${\Bbb S}_{T}$ and ${\Bbb S}_{D}$.

If we put
$$ \sigma_{T\circledcirc D}^{S}=(\mu_{S}\circledcirc T\circledcirc D)\co (S\circledcirc
\psi_{T}^{S}\circledcirc D)\co (\sigma_{T}^{S}\circledcirc
\sigma_{D}^{S})\co (T\circledcirc \tau_{D}^{T}\ot D)$$

we obtain that
\begin{equation}
\label{newsig} \sigma_{T\circledcirc D}^{S}=(\lambda_{1}\circledcirc
\mu_{D})\co (\mu_{T}\circledcirc \lambda_{3}\circledcirc D)\co
(T\circledcirc \lambda_{2}\circledcirc \eta_{S}\circledcirc D)
\end{equation}
and $\sigma_{T\circledcirc D}^{S}$ satisfies (\ref{sigma1}),
(\ref{sigma2}) and (\ref{sigma3}). Indeed: the equality
(\ref{sigma1}) follows by
\begin{equation}
\label{tech2} (\mu_{T}\circledcirc D)\co
(T\circledcirc\lambda_{2})\co (\nabla_{T\circledcirc D}\circledcirc
T)=(\mu_{T}\circledcirc D)\co (T\circledcirc\lambda_{2})
\end{equation}
and the proof for  (\ref{sigma2}) is
\begin{itemize}

\item[ ]$\hspace{0.38cm} \sigma_{T\circledcirc D}^{S}\co
(T\circledcirc D\circledcirc \nabla_{T\circledcirc D}) $

\item [ ]$= (\lambda_{1}\circledcirc D)\co
(\mu_{T}\circledcirc (\lambda_{3}\co
(\mu_{D}\circledcirc\eta_{S})))\co (T\circledcirc
\lambda_{2}\circledcirc D)\co (T\circledcirc D\circledcirc
\nabla_{T\circledcirc D}) $

\item[ ]$=(\lambda_{1}\circledcirc D)\co
(\mu_{T}\circledcirc (\lambda_{3}\co
(\mu_{D}\circledcirc\eta_{S})))\co (T\circledcirc
\lambda_{2}\circledcirc D)$

\item[ ]$=\sigma_{T\circledcirc D}^{S}  $

\end{itemize}
where the first and the third equalities follows by
\begin{equation}
\label{tech3} (S\circledcirc \mu_{D})\co ((\lambda_{3}\co
(D\circledcirc \eta_{S}))\circledcirc D)=\lambda_{3}\co
(\mu_{D}\circledcirc\eta_{S})
\end{equation}
and the second one follows by (\ref{new-nabla}).

Finally, by
$$\nabla_{T\circledcirc D} \co
(T\circledcirc\mu_{D})=(T\circledcirc\mu_{D})\co
(\nabla_{T\circledcirc D}\circledcirc D),$$
$$\nabla_{T\circledcirc D}\co (\mu_{T}\circledcirc D)\co (T\circledcirc \lambda_{2})= 
(\mu_{T}\circledcirc D)\co (T\circledcirc \lambda_{2})$$
and (\ref{YB-Comp}) we obtain (\ref{sigma3}).

As a consequence, by Lemma \ref{iden-ite} and Theorem \ref{prince},
the quadruple
$${\Bbb S}_{T\circledcirc D}=(S, T\circledcirc D,
\psi_{T\circledcirc D}^{S}, \sigma_{T\circledcirc D}^{S}),$$ where
$$\psi_{T\circledcirc D}^{S}=(\lambda_{1}\circledcirc D)\co (T\circledcirc \lambda_{3})\co
(\nabla_{T\circledcirc D}\circledcirc S),$$
 satisfies the equalities (\ref{twis-wcp})
and (\ref{cocy2-wcp}). Then, $(S\circledcirc T\circledcirc  D,
\mu_{S\circledcirc  T\circledcirc D})$ is the iterated weak crossed
product  of $(S\circledcirc  T, \mu_{S\circledcirc T})$ and
$(S\circledcirc D, \mu_{S\circledcirc D})$ with associated
 product
$$\mu_{S\circledcirc T\circledcirc D}=(\mu_S\circledcirc T\circledcirc D)\co (\mu_S\circledcirc
\sigma_{T\circledcirc D}^{S})\co (S\circledcirc \psi_{T\circledcirc
D}^{S}\circledcirc T\circledcirc D),$$ and equivalently
\begin{equation}
\label{product1}\mu_{S\circledcirc T\circledcirc
D}=(\mu_{S}\circledcirc \mu_{T}\circledcirc\mu_{D})\circ
(S\circledcirc ((\lambda_{1}\circledcirc\lambda_{2})\circ
(T\circledcirc\lambda_{3}\circledcirc T)\co (\nabla_{T\circledcirc
D}\circledcirc \nabla_{S\circledcirc T}))\circledcirc D).
\end{equation}

There is no difficulty to extend the previous considerations to the
idempotent clousure $\overline{\mathcal C}$ of an arbitrary
2-category ${\mathcal C}$ and then, the iterated product obtained in
(\ref{product1}) is the one induced by the weak distributive laws in
$\overline{\mathcal C}$ considered by G. B\"{o}hm in
[\cite{Bohm-iterated}, Lemma 2.3].

Also, we have the preunit conditions of Theorem \ref{Teo-preunit}. Indeed, the proof for (\ref{pre-1}) is the following:

\begin{itemize}

\item[ ]$\hspace{0.38cm} (\mu_{S}\circledcirc T\circledcirc D)\co (S\circledcirc \sigma_{T}^{S}\circledcirc D) \co (\lambda_{1}\circledcirc \lambda_{2})\co  (T \circledcirc\lambda_{3}\circledcirc T)\co (\nabla_{T\circledcirc D}\circledcirc \nu_{T})$

\item [ ]$=(S\circledcirc\mu_{T}\circledcirc D)\co (\lambda_{1}\circledcirc\lambda_{2})\co (T\circledcirc \lambda_{3}\circledcirc T)\co  (\nabla_{T\circledcirc D}\circledcirc \nu_{T}) $

\item[ ]$=  (S\circledcirc\mu_{T}\circledcirc D)\co (\lambda_{1}\circledcirc T\circledcirc D)\co (T\circledcirc ((S\circledcirc \lambda_{2})\circ
(\lambda_{3}\circledcirc T)\circ (D\circledcirc
\lambda_{1})))\co (\nabla_{T\circledcirc D}\circledcirc \eta_{T}\circledcirc \eta_{S})  $

\item[ ]$=  (S\circledcirc\mu_{T}\circledcirc D)\co (\lambda_{1}\circledcirc T\circledcirc D)\co (T\circledcirc ((\lambda_{1}\circledcirc D)\circ (T\circledcirc
\lambda_{3})\circ (\lambda_{2}\circledcirc S)))\co (\nabla_{T\circledcirc D}\circledcirc \eta_{T}\circledcirc \eta_{S}) $

\item[ ]$= \nabla_{S\circledcirc T \circledcirc D}\co (\eta_{S}\circledcirc T \circledcirc D) $

\end{itemize}

where the first equality follows by (\ref{idem=idem}) and (\ref{DL3}) for $\lambda_{1}$, the second one follows by the definition of $\nabla_{S\circledcirc T}$, the third one relies on (\ref{YB-Comp}) and the last one is a consequence (\ref{DL1}) for $\lambda_{1}$ and 
$\nabla_{T\circledcirc D}\co \nabla_{T\circledcirc D}=\nabla_{T\circledcirc D}.$

Finally, (\ref{pre-2}) follows by 

\begin{itemize}

\item[ ]$\hspace{0.38cm} (\mu_{S}\circledcirc T\circledcirc D)\co (S \circledcirc \lambda_{1} \circledcirc D)\co (S\circledcirc T\circledcirc \sigma_{D}^{S})\co 
    (S\circledcirc \lambda_{2}\circledcirc D)\co (\nu_{D}\circledcirc T\circledcirc D)$

\item [ ]$= (S\circledcirc T\circledcirc \mu_{D})\co (((S\circledcirc\lambda_{2})\circ
(\lambda_{3}\circledcirc T)\circ (D\circledcirc \lambda_{1}))\circledcirc D)\co (\eta_{D}\circledcirc T\circledcirc \eta_{S}\circledcirc D)$

\item[ ]$=(\lambda_{1}\circledcirc \mu_{D})\co (T\circledcirc\lambda_{3}\circledcirc D)\co (\lambda_{2}\circledcirc S\circledcirc D)\co   (\eta_{D}\circledcirc T\circledcirc \eta_{S}\circledcirc D)  $

\item[ ]$= \nabla_{S\circledcirc T \circledcirc D}\co (\eta_{S}\circledcirc T \circledcirc D) $

\end{itemize}

where  the first equality is a consequence of (\ref{YB-Comp}) and 
$$(\mu_{S} \circledcirc D)\co (S\circledcirc \lambda_{3})\co (\nabla_{S\circledcirc D}\circledcirc S)=(\mu_{S} \circledcirc D)\co (S\circledcirc \lambda_{3}),$$
the second one of (\ref{YB-Comp}) and the last one of (\ref{tech3}).

Therefore, 
$$\nu_{T \circledcirc D}= \nabla_{S\circledcirc T\circledcirc D}\co (\mu_{S}\circledcirc T\circledcirc D)\co (S\circledcirc \lambda_{1}\circledcirc D)\co (\nu_{T}\circledcirc \nu_{ D})$$
is a preunit for $\mu_{S\circledcirc T\circledcirc D}$ and we have that 
$$\nu_{T \circledcirc D}=(\lambda_{1}\circledcirc D)\co (T\circledcirc \lambda_{3})\co (\lambda_{2}\circledcirc S)\co (\eta_{D}\circledcirc \eta_{T}\circledcirc \eta_{S}).$$

}
\end{exemplo}

\begin{exemplo}
\label{ej3} {\rm In this example we will show that the iteration  process proposed recently by D\u{a}u\c{s} and Panaite in  \cite{Pan-2} for Brzezi\'nski's crossed products, is a particular case of the weak iterated products defined in this paper. First we recall from \cite{tb-crpr} the construction of Brzezi\'nski's crossed product in a strict monoidal category: Let $(A,\eta_{A},\mu_{A})$ be a monoid and $V$ an object equipped with a distinguished morphism $\eta_{V}:K\rightarrow V$. Then the object $A\ot V$ is a monoid with unit $\eta_{A}\ot \eta_{V}$ and whose product has the property $\mu_{A\ot V}\co (A\ot \eta_{V}\ot A\ot V)=\mu_{A}\ot V$, if and only if there exists two morphisms $\psi_{V}^{A}:V\ot A\rightarrow A\ot V$, $\sigma_{V}^{A}:V\ot V\rightarrow A\ot V$ satisfying  (\ref{wmeas-wcp}), the twisted condition (\ref{twis-wcp}), the cocycle condition (\ref{cocy2-wcp}) and  
\begin{equation}
\label{brz1} 
\psi_{V}^{A}\co (\eta_{V}\ot A)=A\ot \eta_{V}, 
\end{equation}
\begin{equation}
\label{brz2} 
\psi_{V}^{A}\co (V\ot \eta_{A})=\eta_{A}\ot V, 
\end{equation}
\begin{equation}
\label{brz3} 
\sigma_{V}^{A}\co (\eta_{V}\ot V)=\sigma_{V}^{A}\co (V\ot \eta_{V})=\eta_{A}\ot V.
\end{equation}
If this is the case, the product of $A\ot V$ is the one defined in (\ref{prod-todo-wcp}).  Note that Brzezi\'nski's crossed products are examples of weak crossed products where the associated idempotent is the identity, that is, $\nabla_{A\ot V}=id_{A\ot V}$. Also, in this case the preunit $\nu=\eta_{A}\ot \eta_{V}$ is a unit. 

Given two  Brzezi\'nski's crossed products for  $A\ot V$ and $A\ot W$, in \cite{Pan-2} a new crossed product is defined in $A\ot V\ot W$ if there exists a  morphism $\tau_{W}^{V}:W\ot V\rightarrow V\ot W$ satisfying the condition (i) of Definition \ref{wcp-def} and 
\begin{equation}
\label{DP1} 
(A\ot \tau_{W}^{V})\co (\psi_{W}^{A}\ot V)\co (W\ot\sigma_{V}^{A})=(\sigma_{V}^{A}\ot W)\co (V\ot \tau_{W}^{V})\co (\tau_{W}^{V}\ot V), 
\end{equation}
\begin{equation}
\label{DP2} 
(\psi_{V}^{A}\ot W)\co (V\ot \sigma_{W}^{A})\co (\tau_{W}^{V}\ot W)\co (W\ot \tau_{W}^{V})=(A\ot \tau_{W}^{V})\co (\sigma_{W}^{A}\ot V),
\end{equation}
\begin{equation}
\label{DP3} 
\tau_{W}^{V}\co (\eta_{W}\ot V)=V\ot \eta_{W}, 
\end{equation}
\begin{equation}
\label{DP4} 
\tau_{W}^{V}\co (W\ot \eta_{V})=\eta_{V}\ot W.
\end{equation}

In this case, $\psi_{V\ot W}^{A}=(\psi_{V}^{A}\ot W)\co (V\ot \psi_{W}^{A})$, $\sigma_{V\ot W}^{A}$ is defined as in (\ref{def-sigma}) and  $\eta_{V\ot W}=\eta_{V}\ot \eta_{W}$.  

Under these conditions, the equality (ii) of Definition \ref{wcp-def} holds because:

\begin{itemize}

\item[ ]$\hspace{0.38cm} (\mu_{A}\ot V\ot W)\co (A\ot
\sigma_{V}^{A}\ot W)\co (\psi_{V}^{A}\ot \tau_{W}^{V})\co (V\ot
\sigma_{W}^{A}\ot V)\co (\tau_{W}^{V}\ot W\ot V)$

\item [ ]$=(( (\mu_{A}\ot V)\co (A\ot \sigma_{V}^{A})\co (\psi_{V}^{A}\ot V)\co (V\ot \psi_{V}^{A}))\ot W)\co 
(V\ot V\ot \sigma_{W}^{A})\co (V\ot \tau_{W}^{V}\ot W)\co (\tau_{W}^{V}\ot \tau_{W}^{V}) $

\item[ ]$= (( (\mu_{A}\ot V)\co (A\ot \psi_{V}^{A})\co (\sigma_{V}^{A}\ot A))\ot W)\co 
(V\ot V\ot \sigma_{W}^{A})\co (V\ot \tau_{W}^{V}\ot W)\co (\tau_{W}^{V}\ot \tau_{W}^{V}) $

\item[ ]$=  (\mu_{A}\ot V\ot W)\co (A\ot \psi_{V}^{A}\ot
W)\co (A\ot V\ot \sigma_{W}^{A})\co (A\ot \tau_{W}^{V} \ot W)\co
(\psi_{W}^{A}\ot V\ot W)\co (V\ot \sigma_{V}^{A}\ot W)$
\item[ ]$\hspace{0.38cm}\co (W\ot V\ot
\tau_{W}^{V})$

\end{itemize}

where the first equality follows by (\ref{DP2}), the second one by (\ref{twis-wcp}) for the quadruple ${\Bbb A}_{V}=(A,V,\psi_{V}^{A}, \sigma_{V}^{A})$ and the last one by (\ref{DP1}).

Therefore $\tau_{W}^{V}$ is a twisting morphism and if we consider the link morphism $\Delta_{V\ot W}=id_{V\ot W}$, we obtain that the iterated crossed product proposed in \cite{Pan-2} is a particular instance of the iterated weak crossed product introduced in Theorem \ref{prince}.  Moreover, in this setting,  if the equality (ii) of Definition \ref{wcp-def} holds, composing with $W\ot V\ot \eta_{W}\ot V$ in both sides we obtain  (\ref{DP1}),  and composing with $W\ot \eta_{V}\ot W\ot V$ we obtain (\ref{DP2}). 

}
\end{exemplo}

\section{A different characterization of the iterated weak crossed product}

In this section we obtain a new characterization of the iteration process  following Theorem 1.4 of \cite{mra-proj}. This theorem asserts  the following:

\begin{teo}
\label{uni-1-short} Let $T$ and $B$ be a monoids  in ${\mathcal C}$. Then the
following are equivalent:

\begin{itemize}

\item[(i)] There exist a weak crossed product $(B\otimes W, \mu_{B\otimes W})$ with
preunit $\nu$ and an isomorphism of monoids $\omega:B\times
W\rightarrow T$.
\item[(ii)] There exist an algebra $B$, an object $W$, morphisms $$i_{B}:B\rightarrow
T,\;\;\;\; i_{W}:W\rightarrow T,\;\;\;\;\nabla_{B\otimes W}:
B\otimes W\rightarrow B\otimes W, \;\;\;\;\omega:B\times
W\rightarrow T$$ such that $i_{B}$ is a monoid morphism,
$\nabla_{B\otimes W}$ is an idempotent morphism  of left $B$-modules
for the action $\varphi_{B\otimes W}=\mu_{B}\otimes W$, and $\omega$
is an isomorphism such that
$$\omega\circ p_{B\otimes W}=\mu_{T}\circ (i_{B}\otimes i_{W})$$
 where $B\times W$ is the image of
$\nabla_{B\otimes W}$ and $p_{B\otimes W}$ is the associated
projection.
\end{itemize}
\end{teo}

\begin{teo}
Let ${\Bbb A}_{V}=(A, V, \psi_{V}^{A}, \sigma_{V}^{A})$, ${\Bbb
A}_{W}=(A, W, \psi_{W}^{A}, \sigma_{W}^{A})$ be two quadruples
satisfying  (\ref{twis-wcp}) and (\ref{cocy2-wcp}) with a link
morphism $\Delta_{V\ot W}:V\ot W\rightarrow V\ot W$ and with a
twisting morphism $\tau_{W}^{V}:W\ot V\rightarrow V\ot W$ between
them. Let $(A\ot V, \mu_{A\ot V})$ and $(A\ot W, \mu_{A\ot W})$ be the
weak crossed products associated to ${\Bbb A}_{V}$ and ${\Bbb
A}_{W}$  and suppose that $\nu_{V}:K\rightarrow A\ot V$
and $\nu_{W}:K\rightarrow A\ot W$ are preunits for $\mu_{A\ot V}$
and $\mu_{A\ot W}$. Assume that the morphism  $\sigma_{V\otimes
W}^{A}$, defined in (\ref{def-sigma}), satisfies (\ref{sigma1}),
(\ref{sigma2}), (\ref{sigma3}) and assume also that the  equalities 
(\ref{pre-1}) and (\ref{pre-2}) hold.  

\begin{itemize}

\item[(i)] Let $i_{A\times V}:A\times V\rightarrow A\times (V\ot W)$ be the morphism defined by 
$$i_{A\times V}=p_{A\ot V\ot W}\co (\mu_{A}\ot V\ot W)\co (A\ot \psi_{V}^{A}\ot W)\co (i_{A\ot V}\ot v_{W}),$$
where $A\times (V\ot W)$ is the image of the idempotent morphism $\nabla_{A\ot V\ot W}$ introduced in Definition \ref{link} and $p_{A\ot V\ot W}$  its associated projection. 

If the equality 
\begin{equation}
\label{new-it-1}
\nabla_{A\ot V\ot W}\co (((\mu_{A}\ot V)\co (A\ot \psi_{V}^{A})\co (\sigma_{V}^{A}\ot A))\ot W)\co (V\ot V\ot \nu_{W})
\end{equation}
$$=\nabla_{A\ot V\ot W}\co (((\mu_{A}\ot V)\co (A\ot \sigma_{V}^{A}))\ot W)\co (\psi_{V}^{A}\ot \tau_{W}^{V})\co (V\ot \nu_{W}\ot V),$$
holds, $i_{A\times V}$ is a monoid morphism.

\item[(ii)] If $A\times V$, $p_{A\times V}$ and $i_{A\times V}$ are the image, the projection and the injection associated a $\nabla_{A\otimes V}$, the morphism $\nabla_{(A\times V)\ot W}:(A\times V)\ot W\rightarrow (A\times V)\ot W$ defined by 
$$\nabla_{(A\times V)\ot W}= (p_{A\ot V}\ot W)\co \nabla_{A\ot V\ot W}\co (i_{A\ot V}\ot W),$$
 is idempotent. Moreover, if the following identity holds 
\begin{equation}
\label{new-it-2}
\nabla_{A\ot V\ot W}\co (\sigma_{V}^{A}\ot W)=(((\mu_{A}\ot V)\co (A\ot \psi_{V}^{A}))\ot W)\co (\sigma_{V}^{A}\ot \psi_{W}^{A})\co (V\ot \Delta_{V\ot W}\ot \eta_{A}),
\end{equation}
$\nabla_{(A\times V)\ot W}$ is a morphism of left $A\times V$-modules for $\varphi_{(A\times V)\otimes W}=\mu_{A\times V}\otimes W$.

\item[(iii)] The morphism $\omega :(A\times V)\times W\rightarrow A\times (V\ot W)$ defined by 
$$\omega=p_{A\ot V\ot W}\co (i_{A\ot V}\ot W)\co i_{(A\times V)\ot W},$$
where $i_{(A\times V)\ot W}$ is the injection associated to $\nabla_{(A\times V)\ot W}$, is an isomorphism. Moreover, if the equality 
 \begin{equation}
\label{new-it-3}
\nabla_{A\ot V\ot W}\co (\psi_{V}^{A}\ot W)\co (V\ot \sigma_{W}^{A})=(\psi_{V}^{A}\ot W)\co (V\ot \sigma_{W}^{A})\co (\Delta_{V\ot W}\ot W)
\end{equation}
holds,  then 
$$\omega\circ p_{(A\times V)\otimes W}=\mu_{A\times (V\ot W)}\circ (i_{A\times V}\otimes i_{W})$$
\end{itemize}
for 
$$i_{W}=p_{A\ot V\ot W}\co (\nu_{V}\ot W).$$

Therefore, if (\ref{new-it-1}), (\ref{new-it-2}) and (\ref{new-it-3}) hold, $A\times (V\ot W)$ and $(A\times V)\times W$ are isomorphic as monoids.

\end{teo}

\begin{proof} The proof for (i) is the following:
\begin{itemize}

\item[ ]$\hspace{0.38cm} \mu_{A\times (V\ot W)}\co (i_{A\times V}\ot i_{A\times V}) $

\item[ ]$= p_{A\ot V\ot W}\co (\mu_{A}\ot V\ot W)\co (\mu_{A}\ot \sigma_{V\ot W}^{A})\co (A\ot \psi_{V\ot W}^{A}\ot V\ot W)$
\item[ ]$\hspace{0.38cm} \co (((\mu_{A}\ot V\ot W)\co (A\ot \psi_{V}^{A}\ot W)\co (i_{A\ot V}\ot \nu_{W}))\ot 
((\mu_{A}\ot V\ot W)\co (A\ot \psi_{V}^{A}\ot W)\co (i_{A\ot V}\ot \nu_{W})))$

\item[ ]$= p_{A\ot V\ot W}\co (\mu_{A}\ot V\ot W)\co (\mu_{A}\ot ((\mu_{A}\ot V\ot W)\co (A\ot \psi_{V}^{A}\ot W)\co (\sigma_{V}^{A}\ot \sigma_{W}^{A})))\co (\mu_{A}\ot \psi_{V}^{A} \ot V\ot W\ot W)$
\item[ ]$\hspace{0.38cm}\co (A\ot \psi_{V}^{A}\ot ((A\ot \tau_{W}^{V})\co (\psi_{W}^{A}\ot V)\co (W\ot \psi_{V}^{A}))\ot W)\co (A\ot V\ot \psi_{W}^{A}\ot V\ot A\ot W)\co  (((\mu_{A}\ot \Delta_{V\ot W})$
\item[ ]$\hspace{0.38cm} \co (A\ot \psi_{V}^{A}\ot W)\co (i_{A\ot V}\ot \nu_{W})) \ot i_{A\ot W}\ot \nu_{W})$

\item[ ]$= p_{A\ot V\ot W}\co (\mu_{A}\ot V\ot W)\co (A\ot \mu_{A}\ot V\ot W)\co (A\ot A\ot \psi_{V}^{A}\ot W) $
\item[ ]$\hspace{0.38cm} \co (\mu_{A}\ot ((\mu_A\ot V)\co (A\ot \sigma_{V}^{A})\co (\psi_{V}^{A}\ot V)\co (V\ot \psi_{V}^{A}))\ot \sigma_{W}^{A})\co (A\ot A\ot V\ot V\ot \psi_{W}^{A}\ot W) $
\item[ ]$\hspace{0.38cm} \co (\mu_{A}\ot A\ot V\ot \tau_{W}^{V}\ot \nu_{W}) \co (A\ot A\ot ((\psi_{V}^{A}\ot W)\co (V\ot \psi_{W}^{A})\co (\Delta_{V\ot W}\ot A))\ot V) $
\item[ ]$\hspace{0.38cm}\co  (( (A\ot \psi_{V}^{A}\ot W)\co (i_{A\ot V}\ot \nu_{W}))\ot i_{A\ot V})$

\item[ ]$= p_{A\ot V\ot W}\co (\mu_{A}\ot V\ot W)\co (A\ot \mu_{A}\ot V\ot W)\co (A\ot A\ot \psi_{V}^{A}\ot W) $
\item[ ]$\hspace{0.38cm} \co (\mu_{A}\ot ((\mu_A\ot V)\co (A\ot \psi_{V}^{A})\co (\sigma_{V}^{A}\ot V))\ot \sigma_{W}^{A})\co (A\ot A\ot V\ot V\ot \psi_{W}^{A}\ot W) $
\item[ ]$\hspace{0.38cm} \co (\mu_{A}\ot A\ot V\ot \tau_{W}^{V}\ot \nu_{W}) \co (A\ot A\ot (\nabla_{A\ot V\ot W}\co (\psi_{V}^{A}\ot W)\co (V\ot \psi_{W}^{A}))\ot V) $
\item[ ]$\hspace{0.38cm}\co  (( (A\ot \psi_{V}^{A}\ot W)\co (i_{A\ot V}\ot \nu_{W}))\ot i_{A\ot V})$

\item[ ]$=  p_{A\ot V\ot W}\co (\mu_{A}\ot V\ot W)\co (A\ot \mu_{A}\ot V\ot W)\co (A\ot A\ot \psi_{V}^{A}\ot W)  $
\item[ ]$\hspace{0.38cm} (\mu_{A}\ot \sigma_{V}^{A}\ot ((\mu_{A}\ot W)\co (A\ot \sigma_{W}^{A})\co (\psi_{W}^{A}\ot W)\co (W\ot \nu_{W})))\co (A\ot \mu_{A}\ot V\ot \tau_{W}^{V}) $
\item[ ]$\hspace{0.38cm}\co (A\ot A\ot (\nabla_{A\ot V\ot W}\co (\psi_{V}^{A}\ot W)\co (V\ot \psi_{W}^{A}))\ot V) \co  (( (A\ot \psi_{V}^{A}\ot W)\co (i_{A\ot V}\ot \nu_{W}))\ot i_{A\ot V})$

\item[ ]$=  p_{A\ot V\ot W}\co (\mu_{A}\ot V\ot W)\co (A\ot \mu_{A}\ot V\ot W)\co (A\ot A\ot \psi_{V}^{A}\ot W)  \co (\mu_{A}\ot \sigma_{V}^{A}\ot (\psi_{W}^{A}\co (W\ot \eta_{A})))$
\item[ ]$\hspace{0.38cm}\co (A\ot A\ot V\ot \tau_{W}^{V})\co (A\ot (\nabla_{A\ot V\ot W}\co (A\ot \psi_{V}^{A}\ot W)\co (V\ot ((\mu_{A}\ot W)\co (A\ot \psi_{W}^{A})\co (\nu_{W}\ot A))\ot V)$
\item[ ]$\hspace{0.38cm}\co  ( i_{A\ot V}\ot i_{A\ot V})$

\item[ ]$=  p_{A\ot V\ot W}\co (\mu_{A}\ot V\ot W)\co (A\ot \psi_{V}^{A}\ot W)  \co (A\ot V\ot  (\psi_{W}^{A}\co (W\ot \eta_{A})))$
\item[ ]$\hspace{0.38cm}  \co (\mu_{A}\ot V\ot W)\co (A\ot \sigma_{V}^{A}\ot W)\co (A\ot V\ot \tau_{W}^{V})\co (\nabla_{A\ot V\ot W}\ot V) $
\item[ ]$\hspace{0.38cm}\co (\mu_{A}\ot V\ot W\ot V)\co (A\ot \psi_{V}^{A}\ot W\ot V)\co (i_{A\ot V}\ot ((\beta_{\nu_{W}}\ot V)\co i_{A\ot V}))  $

\item[ ]$=  p_{A\ot V\ot W}\co (\mu_{A}\ot V\ot W)\co (A\ot \psi_{V}^{A}\ot W)  \co (A\ot V\ot  (\psi_{W}^{A}\co (W\ot \eta_{A})))$
\item[ ]$\hspace{0.38cm}  \co (\mu_{A}\ot V\ot W)\co (\mu_{A}\ot ((\mu_{A}\ot V\ot W)\co (A\ot \sigma_{V}^{A}\ot W)\co (\psi_{V}^{A}\ot \tau_{W}^{V})\co (V\ot \sigma_{W}^{A}\ot V)\co (\tau_{W}^{V}\ot W\ot V)))$
\item[ ]$\hspace{0.38cm} \co (A\ot \nu_{W}\ot V\ot W\ot V)\co (\mu_{A}\ot V\ot W\ot V)\co (A\ot \psi_{V}^{A}\ot W\ot V)\co (i_{A\ot V}\ot ((\beta_{\nu_{W}}\ot V)\co i_{A\ot V}))$

\item[ ]$= p_{A\ot V\ot W}\co (\mu_{A}\ot V\ot W)\co (A\ot ((\mu_{A}\ot V)\co (A\ot \psi_{V}^{A}))\ot W)  \co (\mu_{A}\ot A\ot V\ot (\psi_{W}^{A}\co (W\ot \eta_{A})))$
\item[ ]$\hspace{0.38cm}(A\ot ((\mu_{A}\ot ((\psi_{V}^{A}\ot W)\co (V\ot \sigma_{W}^{A})\co (\tau_{W}^{V}\ot W)))\co (A\ot \psi_{W}^{A}\ot V\ot W)\co (\nu_{W}\ot \sigma_{V}^{A}\ot W)\co (V\ot \tau_{W}^{V})))  $
\item[ ]$\hspace{0.38cm}\co (\mu_{A}\ot V\ot W\ot V)\co (A\ot \psi_{V}^{A}\ot W\ot V)\co (i_{A\ot V}\ot ((\beta_{\nu_{W}}\ot V)\co i_{A\ot V})) $

\item[ ]$=p_{A\ot V\ot W}\co (\mu_{A}\ot V\ot W)\co (A\ot \psi_{V}^{A}\ot W)\co (\mu_{A}\ot V\ot (\nabla_{A\ot W}\co \sigma_{W}^{A}))  $
\item[ ]$\hspace{0.38cm}  \co (A\ot ((\mu_{A}\ot \tau_{W}^{V})\co (A\ot \psi_{W}^{A}\ot V)\co (\nu_{W}\ot \sigma_{V}^{A}))\ot W)\co (A\ot V\ot \tau_{W}^{V}) $
\item[ ]$\hspace{0.38cm} \co (\mu_{A}\ot V\ot W\ot V)\co (A\ot \psi_{V}^{A}\ot W\ot V)\co (i_{A\ot V}\ot ((\beta_{\nu_{W}}\ot V)\co i_{A\ot V})) $

\item[ ]$=p_{A\ot V\ot W}\co (\mu_{A}\ot V\ot W)\co (A\ot \psi_{V}^{A}\ot W)\co (\mu_{A}\ot V\ot  \sigma_{W}^{A})  $
\item[ ]$\hspace{0.38cm}  \co (A\ot ((\mu_{A}\ot \tau_{W}^{V})\co (A\ot \beta_{\nu_{W}}\ot V)\co (A\ot \sigma_{V}^{A}))\ot W)\co (A\ot A\ot V\ot \tau_{W}^{V}) $
\item[ ]$\hspace{0.38cm} \co (A\ot \psi_{V}^{A}\ot W\ot V)\co (i_{A\ot V}\ot ((\beta_{\nu_{W}}\ot V)\co i_{A\ot V})) $

\item[ ]$=p_{A\ot V\ot W}\co (\mu_{A}\ot V\ot W)\co (\mu_{A}\ot \mu_{A}\ot V\ot W)  $
\item[ ]$\hspace{0.38cm} \co (A\ot A\ot A\ot ((\mu_{A}\ot V\ot W)\co (A\ot\psi_{V}^{A}\ot W)\co (A\ot V\ot \sigma_{W}^{A})\co (A\ot \tau_{W}^{V}\ot W)\co (\nu_{W}\ot V\ot W))) $
\item[ ]$\hspace{0.38cm} \co (A\ot \mu_{A}\ot \sigma_{V}^{A}\ot W)\co (A\ot A\ot \psi_{V}^{A}\ot \tau_{W}^{V})\co (A\ot \psi_{V}^{A}\ot \nu_{W}\ot V)\co (i_{A\ot V}\ot i_{A\ot V})$

\item[ ]$=p_{A\ot V\ot W}\co (\mu_{A}\ot V\ot W)\co (\mu_{A}\ot \mu_{A}\ot V\ot W)  $
\item[ ]$\hspace{0.38cm} \co (A\ot A\ot A\ot (\psi_{V\ot W}^{A}\co (V\ot W\ot \eta_{A}))) $
\item[ ]$\hspace{0.38cm} \co (A\ot \mu_{A}\ot \sigma_{V}^{A}\ot W)\co (A\ot A\ot \psi_{V}^{A}\ot \tau_{W}^{V})\co (A\ot \psi_{V}^{A}\ot \nu_{W}\ot V)\co (i_{A\ot V}\ot i_{A\ot V})$

\item[ ]$=p_{A\ot V\ot W}\co (\mu_{A}\ot V\ot W) $
\item[ ]$\hspace{0.38cm} (A\ot (\nabla_{A\ot V\ot W}\co (\mu_{A}\ot V\ot W)\co (A\ot \sigma_{V}^{A}\ot W)\co (\psi_{V}^{A}\ot \tau_{W}^{V})\co (V\ot \nu_{W}\ot V)))$
\item[ ]$\hspace{0.38cm} \co (\mu_{A}\ot V\ot V)\co  (A\ot \psi_{V}^{A}\ot V)\co (i_{A\ot V}\ot i_{A\ot V})$

\item[ ]$=p_{A\ot V\ot W}\co (\mu_{A}\ot V\ot W) $
\item[ ]$\hspace{0.38cm} (A\ot (\nabla_{A\ot V\ot W}\co (\mu_{A}\ot V\ot W)\co (A\ot \psi_{V}^{A}\ot W)\co (\sigma_{V}^{A}\ot \nu_{W})))$
\item[ ]$\hspace{0.38cm} \co (\mu_{A}\ot V\ot V)\co  (A\ot \psi_{V}^{A}\ot V)\co (i_{A\ot V}\ot i_{A\ot V})$

\item[ ]$=i_{A\times V}\co \mu_{A\times V}. $

\end{itemize}

The first equality follows because $\mu_{A\ot V\ot W}$ is normalized for  $\nabla_{A\ot V\ot W}$, the second one relies on (\ref{wmeas-wcp}) for ${\Bbb A}_{W}$ and ${\Bbb A}_{V}$, and the third one follows by (i) of Definition \ref{wcp-def} and the associativity of $\mu_{A}$. In the fourth one we applied (\ref{falso-idemp2}) and (\ref{twis-wcp}) for ${\Bbb A}_{V}$. The fifth one follows by (\ref{wmeas-wcp}) for ${\Bbb A}_{V}$ and the associativity of $\mu_{A}$; the sixth one follows by the left linearity for $\nabla_{A\ot V\ot W}$, (\ref{pre1-wcp}) for $\nu_{W}$ and (\ref{wmeas-wcp}) for ${\Bbb A}_{V}$; the seventh one  follows by (\ref{pre3-wcp}) for $\nu_{W}$, the left linearity for $\nabla_{A\ot V\ot W}$ and the associativity of $\mu_{A}$. The eighth one relies on (\ref{pre-2}) and the associativity of $\mu_{A}$, the ninth one is a consequence of 
(ii) of Definition \ref{wcp-def} and the associativity of $\mu_{A}$, and the tenth one follows by (\ref{wmeas-wcp}) for ${\Bbb A}_{V}$. In the eleventh one we used  (\ref{pre3-wcp}) for $\nu_{W}$, (\ref{pre3-wcp}) for $\nu_{W}$ and 
(\ref{idemp-sigma-inv}) for $\sigma_{W}^{A}$. The twelfth  one follows by  (\ref{wmeas-wcp}) for ${\Bbb A}_{V}$ and the associativity of $\mu_{A}$; the thirteenth one follows by (\ref{pre-2}) and the fourteenth one follows by the left linearity for $\nabla_{A\ot V\ot W}$. The fifteenth one is a consequence of (\ref{new-it-1}) and the last one follows by the associativity of $\mu_{A}$, the left linearity for $\nabla_{A\ot V\ot W}$ and by (\ref{fi-nab}) for ${\Bbb A}_{V}$.

Therefore, $i_{A\times V}$ is multiplicative and, by (\ref{preunit-idemp}), we have  
$$i_{A\times V}\co \eta_{A\times V}=p_{A\ot V\ot W}\co (\mu_{A}\ot V\ot W)\co (A\ot \psi_{V}^{A}\ot W)\co ((\nabla_{A\ot V}\co \nu_{V})\ot \nu_{W})$$
$$=p_{A\ot V\ot W}\co (\mu_{A}\ot V\ot W)\co (A\ot \psi_{V}^{A}\ot W)\co ( \nu_{V}\ot \nu_{W})=\eta_{A\times (V\ot W)}.$$

(ii) The morphism $\nabla_{(A\times V)\ot W}= (p_{A\ot V}\ot W)\co \nabla_{A\ot V\ot W}\co (i_{A\ot V}\ot W),$ is idempotent because 

\begin{itemize}

\item[ ]$\hspace{0.38cm}  \nabla_{(A\times V)\ot W}\co \nabla_{(A\times V)\ot W}$

\item[ ]$= (p_{A\ot V}\ot W)\co \nabla_{A\ot V\ot W}\co ((\nabla_{A\ot V}\co (\mu_{A}\ot V)\co (A\ot \psi_{V}^{A}))\ot W)\co (A\ot V\ot (\psi_{W}^A\co (W\ot \eta_{A})))$
\item[ ]$\hspace{0.38cm}\co (A\ot \Delta_{V\ot W})\co (i_{A\ot V}\ot W) $

\item[ ]$= (p_{A\ot V}\ot W)\co (\mu_{A}\ot V\ot W)\co (A\ot \psi_{V}^{A}\ot W)\co (\mu_{A}\ot V\ot (\psi_{W}^A\co (W\ot \eta_{A})))\co (A\ot ((A\ot \Delta_{V\ot W})$
\item[ ]$\hspace{0.38cm}\co (\psi_{V\ot W}^{A}\co (V\ot W\ot \eta_{A}))))\co (i_{A\ot V}\ot W) $

\item[ ]$=  (p_{A\ot V}\ot W)\co (\mu_{A}\ot V\ot W)\co (A\ot ((\mu_{A}\ot V\ot W)\co (A\ot \psi_{V}^{A}\ot W)\co (\psi_{V}^{A}\ot \psi_{W}^{A})\co (V\ot \psi_{W}^{A} \ot A)$
\item[ ]$\hspace{0.38cm}\co (V\ot W\ot \eta_{A}\ot \eta_{A})))\co (A\ot \Delta_{V\ot W})\co (i_{A\ot V}\ot W) $

\item[ ]$=\nabla_{(A\times V)\ot W},$

\end{itemize}

where the first equality follows by definition, the second one follows by (\ref{fi-nab}), the third one relies on (\ref{falso-idemp}), and the las one follows by  (\ref{wmeas-wcp}) for ${\Bbb A}_{V}$ and ${\Bbb A}_{W}$.

On the other hand, 

\begin{itemize}

\item[ ]$\hspace{0.38cm}  \nabla_{(A\times V)\ot W}\co \varphi_{(A\times V)\ot W}$

\item[ ]$= (p_{A\ot V}\ot W)\co  \nabla_{A\ot V\ot W}\co ((\nabla_{A\ot V}\co (\mu_{A}\ot V)\co (\mu_{A}\ot \sigma_{V}^{A})\co (A\ot \psi_{V}^{A}\ot V)\co (i_{A\ot V}\ot i_{A\ot V}))\ot W)$

\item[ ]$= (p_{A\ot V}\ot W)\co (\mu_{A}\ot V\ot W)\co (\mu_{A}\ot (\nabla_{A\ot V\ot W}\co (\sigma_{V}^A\ot W)))\co  (((A\ot \psi_{V}^{A}\ot V)\co (i_{A\ot V}\ot i_{A\ot V}))\ot W)$

\item[ ]$=(p_{A\ot V}\ot W)\co (\mu_{A}\ot V\ot W)\co (\mu_{A}\ot ((\mu_{A}\ot V\ot W)\co (A\ot \psi_{V}^{A}\ot W)\co (\sigma_{V}^{A}\ot (\psi_{W}^{A}\co (W\ot \eta_{A})))$
\item[ ]$\hspace{0.38cm} \co (V\ot \Delta_{V\ot W})))\co  (((A\ot \psi_{V}^{A}\ot V)\co (i_{A\ot V}\ot i_{A\ot V}))\ot W)$

\item[ ]$= (p_{A\ot V}\ot W)\co (\mu_{A}\ot V\ot W)\co (\mu_{A}\ot ((\mu_{A}\ot V)\co (A\ot \sigma_{V}^{A})\co 
(\psi_{V}^{A}\ot V)\co (V\ot \psi_{V}^{A}))\ot W)  $
\item[ ]$\hspace{0.38cm} \co (A\ot \psi_{V}^{A}\ot V\ot (\psi_{W}^{A}\co (W\ot \eta_{A})))\co (A\ot V\ot A\ot \Delta_{V\ot W})\co (i_{A\ot V}\ot i_{A\ot V}\ot W)$

\item[ ]$=\varphi_{(A\times V)\ot W}\co (A\times V \ot \nabla_{(A\times V)\ot W})$

\end{itemize}
where the first equality follows by definition, the second one follows by the left linearity of $\nabla_{A\ot V}$ and 
(\ref{idemp-sigma-inv}), the third one relies on  (\ref{new-it-2}) and the fifth one is a consequence of (\ref{fi-nab}), (\ref{wmeas-wcp}) for ${\Bbb A}_{V}$ and the associativity of $\mu_{A}$.

Finally, we will prove (iii). The morphism $\omega=p_{A\ot V\ot W}\co (i_{A\ot V}\ot W)\co i_{(A\times V)\ot W}$ is an isomorphism with inverse 
$$\omega^{-1}=p_{(A\times V)\ot W}\co (p_{A\ot V}\ot W)\co i_{A\ot V\ot W}$$ because 
$$\omega^{-1}\co \omega=p_{(A\times V)\ot W}\co \nabla_{(A\times V)\ot W}\co i_{(A\times V)\ot W}=id_{(A\times V)\ot W}$$
and, by (\ref{fi-nab}), (\ref{falso-idemp2}) and the left linearity of $\nabla_{A\ot V\ot W}$, we have 
\begin{itemize}

\item[ ]$\hspace{0.38cm}  \omega\co \omega^{-1}$

\item[ ]$=p_{A\ot V\ot W}\co ((\nabla_{A\ot V}\co (\mu_{A}\ot V)\co (A\ot \psi_{V}^{A}))\ot W)\co (A\ot V\ot \psi_{W}^{A})\co (A\ot \Delta_{V \ot W}\ot \eta_{A})\co (\nabla_{A\ot V}\ot W)$
\item[ ]$\hspace{0.38cm}\co i_{A\ot V\ot W}$

\item[ ]$=p_{A\ot V\ot W}\co ( ((\mu_{A}\ot V)\co (A\ot \psi_{V}^{A}))\ot W)\co (A\ot V\ot \psi_{W}^{A})\co (A\ot \Delta_{V \ot W}\ot \eta_{A})\co (\nabla_{A\ot V}\ot W)$
\item[ ]$\hspace{0.38cm}\co i_{A\ot V\ot W}$

\item[ ]$=p_{A\ot V\ot W}\co (\mu_{A}\ot V\ot W)\co (A\ot \nabla_{A\ot V\ot W})\co (A\ot \psi_{V}^{A}\ot W)\co (A\ot V\ot ((\psi_{W}^{A}\co (W\ot \eta_{A}))))\co (\nabla_{A\ot V}\ot W) $
\item[ ]$\hspace{0.38cm}\co i_{A\ot V\ot W}$

\item[ ]$=p_{A\ot V\ot W}\co (\mu_{A}\ot V\ot W)\co (A\ot \psi_{V}^{A}\ot W)\co (A\ot V\ot ((\psi_{W}^{A}\co (W\ot \eta_{A}))))\co (\nabla_{A\ot V}\ot W) \co i_{A\ot V\ot W}$

\item[ ]$=p_{A\ot V\ot W}\co (\mu_{A}\ot V\ot W)\co (A\ot \psi_{V}^{A}\ot W)\co (A\ot V\ot ((\psi_{W}^{A}\co (W\ot \eta_{A}))))\co  i_{A\ot V\ot W}$

\item[ ]$=p_{A\ot V\ot W}\co (\mu_{A}\ot V\ot W)\co (A\ot \nabla_{A\ot V\ot W})\co (A\ot \psi_{V}^{A}\ot W)\co (A\ot V\ot ((\psi_{W}^{A}\co (W\ot \eta_{A}))))\co  i_{A\ot V\ot W}$ 

\item[ ]$=p_{A\ot V\ot W}\co (\mu_{A}\ot V\ot W)\co  (A\ot \psi_{V}^{A}\ot W)\co (A\ot V\ot ((\psi_{W}^{A}\co (W\ot \eta_{A}))))\co (A\ot \Delta_{V\ot W})\co  i_{A\ot V\ot W}$

\item[ ]$= id_{A\times (V\ot W)}.$
\end{itemize}

Moreover, if (\ref{new-it-3}) holds, we have the following:

\begin{itemize}

\item[ ]$\hspace{0.38cm}\mu_{A\times (V\ot W)}\circ (i_{A\times V}\otimes i_{W})$

\item[ ]$= p_{A\ot V\ot W}\co (\mu_{A}\ot V\ot W)\co (\mu_{A}\ot \psi_{V}^{A}\ot W)\co (A\ot A\ot V\ot \sigma_{W}^{A}) $
\item[ ]$\hspace{0.38cm}\co (\mu_{A}\ot ((\mu_{A}\ot V\ot W)\co (A\ot \sigma_{V}^{A}\ot W)\co
(\psi_{V}^{A}\ot \tau_{W}^{V})\co (V\ot \psi_{W}^{A}\ot V)\co
(\Delta_{V\ot W}\ot \nu_{V}))\ot W)$
\item[ ]$\hspace{0.38cm}\co (A\ot \psi_{V}^{A}\ot W\ot W)\co (i_{A\ot V}\ot \nu_{W}\ot W)$

\item[ ]$=p_{A\ot V\ot W}\co (\mu_{A}\ot V\ot W)\co (A\ot \mu_{A}\ot V\ot W)\co (A\ot A\ot \psi_{V}^{A}\ot W)\co (\mu_{A}\ot \psi_{V}^{A}\ot \sigma_{W}^{A}) $
\item[ ]$\hspace{0.38cm}\co (A\ot A\ot V\ot (\psi_{W}^{A}\co (W\ot \eta_{A}))\ot W)\co (A\ot A\ot \Delta_{V\ot W}\ot W)\co (A\ot \psi_{V}^{A}\ot W\ot W)\co  (i_{A\ot V}\ot \nu_{W}\ot W)$

\item[ ]$= p_{A\ot V\ot W}\co (\mu_{A}\ot V\ot W)\co (A\ot  \psi_{V}^{A}\ot W)\co (\mu_{A}\ot V\ot ((\mu_{A}\ot W)\co (A\ot \sigma_{W}^{A})\co ( (\psi_{W}^{A}\co (W\ot \eta_{A}))\ot W)$
\item[ ]$\hspace{0.38cm}\co  (A\ot A\ot \Delta_{V\ot W}\ot W)\co (A\ot \psi_{V}^{A}\ot W\ot W)\co  (i_{A\ot V}\ot \nu_{W}\ot W)$

\item[ ]$= p_{A\ot V\ot W}\co (\mu_{A}\ot V\ot W)\co (A\ot  \psi_{V}^{A}\ot W)\co (\mu_{A}\ot V\ot ((\mu_{A}\ot W)\co (A\ot \sigma_{W}^{A})\co (\nabla_{A\ot W}\co (\eta_{A}\ot W))))$
\item[ ]$\hspace{0.38cm}\co  (A\ot A\ot \Delta_{V\ot W}\ot W)\co (A\ot \psi_{V}^{A}\ot W\ot W)\co  (i_{A\ot V}\ot \nu_{W}\ot W)$ 

\item[ ]$=p_{A\ot V\ot W}\co (\mu_{A}\ot V\ot W)\co (\mu_{A}\ot ((\psi_{V}^{A}\ot W)\co (V\ot \sigma_{W}^{A})\co (\Delta_{V\ot W}\ot W))) \co (A\ot \psi_{V}^{A}\ot W\ot W)$
\item[ ]$\hspace{0.38cm} \co  (i_{A\ot V}\ot \nu_{W}\ot W)$ 

\item[ ]$=p_{A\ot V\ot W}\co (\mu_{A}\ot V\ot W)\co (\mu_{A}\ot (\nabla_{A\ot V\ot W}\co (\psi_{V}^{A}\ot W)\co (V\ot \sigma_{W}^{A})))\co (A\ot \psi_{V}^{A}\ot W\ot W) $
\item[ ]$\hspace{0.38cm} \co  (i_{A\ot V}\ot \nu_{W}\ot W)$

\item[ ]$=p_{A\ot V\ot W}\co (\mu_{A}\ot V\ot W)\co  (A\ot \psi_{V}^{A}\ot W)\co (i_{A\ot V}\ot ((\mu_{A}\ot W)\co (A\ot \sigma_{W}^{A})\co (\nu_{W}\ot W))) $

\item[ ]$=p_{A\ot V\ot W}\co (\mu_{A}\ot V\ot W)\co (A\ot (\nabla_{A\ot V\ot W}\co (\psi_{V}^{A}\ot W)\co (V\ot  (\psi_{W}^{A}\co (W\ot \eta_{A})))))\co (i_{A\ot V}\ot W)$

\item[ ]$= p_{A\ot V\ot W}\co (\nabla_{A\ot V}\ot W)\co \nabla_{A\ot V\ot W}\co (i_{A\ot V}\ot W)$ 

\item[ ]$=\omega\circ p_{(A\times V)\otimes W},$
\end{itemize}

where the first equality follows because $\mu_{A\ot V\ot W}$ is normalized for  $\nabla_{A\ot V\ot W}$ and by the associativity of $\mu_{A}$, the second one follows by (\ref{pre-2}) and by the associativity of $\mu_{A}$, the third one relies on (\ref{wmeas-wcp}) for ${\Bbb A}_{V}$ and the fourth one is a consequence of the properties of $\nabla_{A\ot W}$. The fifth one follows by (\ref{aw1}) and by the associativity of $\mu_{A}$, the sixth one follows by (\ref{new-it-3}), and in the seventh one we used the left linearity of $\nabla_{A\ot V}$ and (\ref{wmeas-wcp}) for ${\Bbb A}_{V}$. In the eighth one we applied the left linearity of $\nabla_{A\ot V}$ and (\ref{pre2-wcp}) for $\nu_{W}$. The ninth one follows by (\ref{falso-idemp2}) and by (\ref{fi-nab}), and the last one follows by definition. 

The final assertion of this theorem follows by Theorem \ref{uni-1-short}.

\end{proof}

\begin{exemplo}
{\rm In this example we will see that the equalities (\ref{new-it-1}), (\ref{new-it-2}) and (\ref{new-it-3}) hold in the examples (\ref{ej1}), (\ref{ej2}) and (\ref{ej3}) of the previous section.

For the Example (\ref{ej1}) the identities (\ref{new-it-1}), (\ref{new-it-2}) and (\ref{new-it-3}) hold  because 
$$\psi_{T}^{S}=\lambda_{1},\;\sigma_{T}^{S}= \mu_{T}\circledcirc \eta_{S},\; \tau_{T}^{D}=\lambda_{2},\; \nu_{D}=\eta_{S}\circledcirc\eta_{D},$$
and 
$$\Delta_{T\circledcirc D}=id_{T\circledcirc D},\; \nabla_{S\circledcirc T \circledcirc D}=id_{S\circledcirc T \circledcirc D}.$$

In the case of the Example (\ref{ej2}) we have that 
$$\psi_{T}^{S}=\lambda_{1},\;\sigma_{T}^{S}= \lambda_{1}\co (\mu_{T}\circledcirc \eta_{S}),\; \tau_{T}^{D}=\lambda_{2},\; \nu_{D}=\nabla_{S\circledcirc D}\co (\eta_{S}\circledcirc\eta_{D}),$$
and $\Delta_{T\circledcirc D}=\nabla_{T\circledcirc D}.$ Therefore, by the usual arguments, we obtain that (\ref{new-it-1}), (\ref{new-it-2}) and (\ref{new-it-3}) hold  because  
$$ \nabla_{S\circledcirc T \circledcirc D} \co (\mu_{S}\circledcirc T\circledcirc D)\co (S\circledcirc\sigma_{T}^{S}\circledcirc D)\co (\psi_{S}^{T}\circledcirc \tau_{T}^{D})\co (T\circledcirc\nu_{D}\circledcirc T)$$
$$=(\lambda_{1}\circledcirc D)\co (T\circledcirc \lambda_{3})\co ((\nabla_{T\circledcirc D}\co (\mu_{T}\circledcirc \eta_{D}))\circledcirc \eta_{S})$$
$$=\nabla_{S\circledcirc T \circledcirc D} \co (\mu_{S}\circledcirc T\circledcirc D)\co (S\circledcirc \psi_{T}^{S}\circledcirc D)\co (\sigma_{T}^S\circledcirc \nu_{D}),$$
$$\nabla_{S\circledcirc T \circledcirc D} \co (\sigma_{T}^S\circledcirc D)
=(S\circledcirc \mu_{T}\circledcirc D)\co (\lambda_{1}\circledcirc\lambda_{2})\co (T\circledcirc \lambda_{3}\circledcirc T)\co (\mu_{T}\circledcirc D\circledcirc (\lambda_{1}\co (\eta_{T}\circledcirc \eta_{S})))$$
$$=(\mu_{S}\circledcirc T\circledcirc D)\co (S\circledcirc \psi_{T}^{S}\circledcirc D)\co (\sigma_{T}^{S}\circledcirc \psi_{D}^{S})\co (T\circledcirc \Delta_{T\circledcirc D}\circledcirc \eta_{S}), $$
and 
$$ (\psi_{T}^{S}\circledcirc D)\co (T\circledcirc \sigma_{D}^{S})\co (\Delta_{T\circledcirc D}\circledcirc D)=
(S\circledcirc \mu_{T}\circledcirc D)\co (\lambda_{1}\circledcirc\lambda_{2})\co (T\circledcirc \lambda_{3}\circledcirc T)\co (T\circledcirc \mu_{D}\circledcirc (\lambda_{1}\co (\eta_{T}\circledcirc \eta_{S})))$$
$$ =\nabla_{S\circledcirc T \circledcirc D} \co (\psi_{T}^{S}\circledcirc D)\co (T\circledcirc \sigma_{D}^{S}).$$

Finally, in Example (\ref{ej3}) we have that 
$$\nu_{W}= \eta_{A}\ot \eta_{W},\; \Delta_{V\ot W}=id_{V\ot W},\; \nabla_{A\ot V\ot W}=id_{A\ot V\ot W}, $$
and then (\ref{new-it-1}), (\ref{new-it-2}) and (\ref{new-it-3})  follow easily.

}
\end{exemplo}

\section*{Acknowledgements}
The  first and second authors were supported by  Ministerio de Econom\'{\i}a y Competitividad of Spain (European Feder support included). Grant MTM2013-43687-P: Homolog\'{\i}a, homotop\'{\i}a e invariantes categ\'oricos en grupos y \'algebras no asociativas.


\begin{thebibliography}{99}

\bibitem{AM2} A. Agore, G. Militaru,
{\em Unified products and split extensions of Hopf algebras}.
Contemporary Math. 585 (2013), 1-15.

\bibitem{nmra4} J.N. Alonso Álvarez, J.M. Fernández Vilaboa, R.
González Rodríguez, A.B. Rodríguez Raposo, A.B.,  {\em Crosssed
products in weak contexts}. Appl. Cat. Structures 3 (2010), 231-258.

\bibitem{Beck}
J. Beck, {\em Distributive laws}, Springer LNM 80,
119-140 (1969).

\bibitem{bohm} G. B\"ohm, {\em The weak theory of monads},
 Adv. Math. 225 (2010), 1-32.

\bibitem{Bohm-iterated}
G. B\"ohm. {\em On the iteration of weak wreath products}. Theory
and Applications of Categories 26 (2012), 30-59.

\bibitem{Gabi-Pepe1}
G. B\"ohm, J. G\'omez-Torrecillas, {\em Bilinear factorizations of algebras}. Bulletin of the Belgian Mathematical Society Simon Stevin 20 (2013), 1-24.

\bibitem{Gabi-Pepe2}
G. B\"ohm, J. G\'omez-Torrecillas, {\em On the double crossed product of weak Hopf algebras}. Contemporary Math. 585 (2013) 153-173.

\bibitem{tb-crpr} T. Brzezi\'nski, {\em Crossed products by a
coalgebra}. Comm. in Algebra 25 (1997), 3551-3575.

\bibitem{Caen}
S. Caenepeel, E. De Groot, {\em Modules over weak entwining
structures}. Contemporary Math. 267 (2000), 4701-4735.

\bibitem{CSV} A. Cap, H. Schichl, J. Vanzura, {\em On
twisted tensor products of algebras}. Comm. in Algebra  23 (1995),
4701-4735.

\bibitem{Pan-2}  L. D\u{a}u\c{s},  F. Panaite, {\em A new way to iterate 
Brzezi\'nski crossed products}, arXiv:1502.00031 (2015).


\bibitem{mra-preunit}  J.M. Fern\'andez Vilaboa,
R. Gonz\'alez Rodr\'{\i}guez, A.B. Rodr\'{\i}guez Raposo, {\em
Preunits and weak crossed products}. J. of Pure Appl. Algebra 213
(2009), 2244-2261.

\bibitem{mra-proj} J.M. Fern\'andez Vilaboa, R. Gonz\'alez Rodr\'{\i}guez, A.B.
Rodr\'{\i}guez Raposo,
{\em Weak crossed biproducts and weak projections}. Sci. China Math.
 55  (2012), 1443-11460 

\bibitem{mra-partial-unif} J.M. Fern\'andez Vilaboa, R. Gonz\'alez Rodr\'{\i}guez,
 A.B. Rodr\'{\i}guez Raposo,
{\em Partial and unified crossed products are weak crossed
products}. Contemporary Math. 585 (2013) 261-274.

\bibitem{Pascual-Javi2}
P. Jara Mart\'{\i}nez, J. L\'opez Peña, F. Panaite, F. van
Oystaeyen, {\em On iterated tensor product of algebras}. Internat.
J. Math. 19 (2008), 1053-1101.

\bibitem{Christian} C. Kassel, \emph{Quantum Groups}.  Springer-Verlag,  1995.

\bibitem{LS} S. Lack, R. Street, {\em The formal theory of monads II}.
 J. of Pure Appl. Algebra  175 (2002), 243-265.

\bibitem{Pascual-Javi1}
 J. L\'opez Peña, F. Panaite, F. van
Oystaeyen, {\em General twisting of algebras}. Adv. Math. 212
(2007), 315-337.

\bibitem{Macl}  S. Maclane,  {\em Categories for the
Working Mathematicien}.  Springer-Verlag, 1971.

\bibitem{Maj} S. Majid. {\em More examples of bicrossproduct and double cross product Hopf algebras.} Israel J. Math., 72 (1990), 133-148.

\bibitem{partial}   M. Muniz S. Alves, E. Batista, M. Dokuchaev, A.
Paques,  {\em Twisted partial actions of Hopf algebras},
Israel J. of Math. 197 (2013), 263-308

\bibitem{Pan-1} F. Panaite, {\em Iterated crossed products}. J. Algebra Appl. 13 (2014), 
1450036 (14 pages)

\bibitem{ana1}  A.B. Rodr\'{\i}guez Raposo,  {\em Crossed products for
weak Hopf algebras}. Comm. in Algebra 37 (2009),
2274-2289.

\bibitem{Street-FTM}
R. Street, {\em The formal theory of monads}. J. Pure Appl. Algebra
2 (1972), 149-168.

\bibitem{Street-WDL}
R. Street, {\em Weak distributive laws}. Theory and Applications of
Categories 22 (2009), 313-320.

\bibitem{TAM} D. Tambara,  {\em The endomorphism
bialgebra of an algebra}. J. Fac. Univ. Tokyo Sect. IA  37 (1990),
425-456.


\end{thebibliography}
\end{document}